\input amstex
\documentstyle{amsppt}
\pagewidth{360pt}
\pageheight{583pt}
\nopagenumbers
\accentedsymbol\tx{\tilde x}
\accentedsymbol\ty{\tilde y}
\def\twosum{\operatornamewithlimits{\sum\!\sum}}
\def\const{\operatorname{const}}
\def\End{\operatorname{End}}
\def\Hom{\operatorname{Hom}}
\def\Ker{\operatorname{Ker}}
\def\Cl{\operatorname{Cl}}
\def\tr{\operatorname{tr}}

\def\id{\operatorname{\bold{id}}}
\def\compos{\,\raise 1pt\hbox{$\sssize\circ$} \,}
\rightheadtext{On a point symmetry analysis \dots}
\leftheadtext{V\.~Dmitrieva, E\.~Neufeld, R\.~Sharipov,
A\.~Tsaregorodtsev}
\topmatter
\title On a point symmetry analysis for generalized 
diffusion type equations.
\endtitle
\author
V~\.V\.~Dmitrieva, E\.~G\.~Neufeld, R\.~A\.~Sharipov,
A\.~A\.~Tsaregorodtsev
\endauthor
\address
Department of Mathematics, Bashkir State University,\newline
Frunze street 32, 450074, Ufa, Russia
\endaddress
\email \vtop{\hsize 6cm \noindent DmitrievaVV\@ic.bashedu.ru\newline
R\_\hskip 0.7pt Sharipov\@ic.bashedu.ru\vskip 3pt}
\endemail
\urladdr http://www.geocities.com/CapeCanveral/Lab/5341
\endurladdr
\address
Department of Mathematics, Ufa State Aviation Technical 
University,\newline
Karl Marks street 12, 450000 Ufa, Russia.
\endaddress
\email
anton\@gov.math.ugatu.ac.ru
\endemail
\thanks Authors are grateful to International Soros Foundation 
(``Open Society'' institute), Russian Fund for Fundamental 
Research, and the Academy of Sciences of the Republic Bashkortostan
for financial support in 1998.
\endthanks
\abstract Generalized diffusion type equations are considered and
point symmetry analysis is applied to them. The equations with 
extremal order point symmetry algebras are described. Some old 
geometrical results are rederived in connection with theory of 
these equation.
\endabstract
\endtopmatter
\loadbold
\document
\head
1. Introduction.
\endhead
    In this paper we consider a class of systems of evolutional
equations, which we call {\it the equations of generalized diffusion}:
$$
\frac{\partial y^i}{\partial\,t}=\sum^n_{j=1}A^i_j\left(
\frac{\partial^2 y^j}{\partial x^2}+\shave{\sum^n_{r=1}\sum^n_{s=1}}
\Gamma^j_{rs}\,\frac{\partial y^r}{\partial x}\,\frac{\partial y^s}
{\partial x}\right),\quad i=1,\,\ldots,\,n.\hskip -2em
\tag1.1
$$
Here $x$ and $t$ are independent variables, $y^1,\,\ldots,\,y^n$ are
dependent variables, and coefficients  $A^i_j$ and $\Gamma^j_{rs}$
are the functions of $y^1,\,\ldots,\,y^n$. They do not depend on
$x$ and $t$ explicitly.\par
     Special case of \thetag{1.1} is formed by {\it the equations 
of diffusional type}, which describe the diffusion phenomena in 
multicomponent mixtures:
$$
\frac{\partial y^i}{\partial\,t}=\sum^n_{j=1}\frac{\partial}
{\partial x}\left(A^i_j\,\frac{\partial y^j}{\partial x}
\right),\quad i=1,\,\ldots,\,n,\hskip -2em
\tag1.2
$$
(see paper \cite{1} and reference list there). The equations of magnet 
in Heisenberg model (see \cite{2} and \cite{3}) are also reduced to 
\thetag{1.1} in classical limit:
$$
\bold S_t=[\bold S,\,\bold S_{xx}].
\tag{1.3}
$$
Here $\bold S$ is a unit length three-dimensional vector, which describes
the magnetization of the media. Square brackets denote the vector product.
The equations \thetag{1.3} have multidimensional generalization, where
$\bold S$ is an element of some Lie algebra and square brackets denote 
the commutator in this algebra (see for instance \cite{4} and \cite{5}).
They can be restricted to any orbit of coadjoint action of corresponding
Lie group $G$ on $L$. When written in local coordinates $y^1,\,\ldots,\,
y^n$ on such orbit, the equations \thetag{1.3} have the form \thetag{1.1}
\par
     An important feature of the class of equations \thetag{1.1} is that 
it is closed with respect to the point transformations of the form
$$
\align
\ty^1&=\ty^1(y^1,\ldots,y^n),\hskip -2em\\
.\ .&\ .\ .\ .\ .\ .\ .\ .\ .\ .\ .\ .\ .\ \hskip -2em
\tag1.4\\
\ty^n&=\ty^1(y^1,\ldots,y^n),\hskip -2em
\endalign
$$
which do not change independent variables $x$ and $t$. Therefore we can
treat variables $y^1,\,\ldots,\,y^n$ as local coordinates on some manifold
\footnote{for equations \thetag{1.3} such manifold comes from initial
statement of the problem --- this is the orbit of coadjoint action of
Lie group on its Lie algebra.}. Let's denote it by $M$. Change of
variables \thetag{1.4} corresponds to the transition from one local
map on $M$ to another. By means of direct calculations one derives the
following transformation rules for $A^i_j$ and $\Gamma^j_{rs}$ 
under the point transformation of the form \thetag{1.4}:
$$
\align
&A^k_i=\sum^n_{m=1}\sum^n_{p=1}S^k_m\,T^p_i\,\tilde A^m_p,
\hskip -2em
\tag1.5\\
&\Gamma^k_{ij}=\sum^n_{m=1}\sum^n_{p=1}\sum^n_{q=1}S^k_m\,T^p_i\,
T^q_j\,\tilde\Gamma^m_{pq}+\sum^n_{m=1}S^k_m\,\frac{\partial T^m_i}
{\partial y^j}.\hskip -2em
\tag1.6
\endalign
$$
Here $T$ and $S$ are Jacobi matrices for direct and inverse transformation
\thetag{1.4}:
$$
\xalignat 2
&S^i_j=\frac{\partial y^i}{\partial\ty^j},
&&T^i_j=\frac{\partial\ty^i}{\partial y^j}.
\endxalignat
$$
In geometry it is customary to call $S$ {\it the matrix of direct 
transition}, while $T$ is called {\it the inverse transition 
matrix} (see for instance \cite{6}).\par
     From \thetag{1.5} and \thetag{1.6} one can understand that parameters
$\Gamma^j_{rs}$ define symmetric affine connection $\Gamma$ on $M$ and
parameters $A^i_j$ define a tensor field $\bold A$ of the type $(1,1)$.
This provides us with ample opportunity to apply powerful differential
geometric methods to the study of equations \thetag{1.1}, e\.~g\.
in \cite{1} we have found an effective criterion for checking whether the
equations \thetag{1.1} can be brought to the form \thetag{1.2} by means
of some transformation \thetag{1.4} or not. In this paper we consider
the problem of describing the equations \thetag{1.1} whose point symmetry
algebras are large enough (extremal in some sense).
\head
2. Vector fields and point symmetries.
\endhead
     Let $\boldsymbol\eta$ be a vector field on the manifold $M$. In local
coordinates it is represented by differential operator
$$
\boldsymbol\eta=\sum^n_{i=1}\eta^i\,\frac{\partial}{\partial y^i}.
\hskip -2em
\tag2.1
$$
Vector field \thetag{2.1} gives rise to the local one-parametric group
of transformations $\varphi_\tau\!:M\to M$, while each solution of
\thetag{1.1} is interpreted as two-parametric set of points in $M$
($x$ and $t$ are parameters). Applying $\varphi_\tau$ to this set we
obtain another two-parametric set of points with the same parameters
$x$ and $t$. Vector field $\boldsymbol\eta$ is called {\it a point 
symmetry} for some particular system of equations, if for any solution 
of this system and for arbitrary value of $\tau$ the transformed 
two-parametric set of points satisfies the same equations \thetag{1.1} 
as the initial set. From this statement one can derive the following 
determining equations for the field of point symmetry:
$$
\gather
\sum^n_{k=1}\eta^k\,\frac{\partial A^i_j}{\partial y^k}-\sum^n_{k=1}
A^k_j\,\frac{\partial\eta^i}{\partial y^k}+\sum^n_{k=1}A^i_k\,
\frac{\partial\eta^k}{\partial y^j}=0,\hskip-2em
\tag2.2\\
\vspace{2ex}
\gathered
\sum^n_{j=1}\sum^n_{k=1}\left(\frac{\partial\eta^i}{\partial y^k}\,A^k_j
-\frac{\partial A^i_j}{\partial y^k}\,\eta^k\right)\Gamma^j_{rs}=\\
=\sum^n_{j=1}\sum^n_{k=1}A^i_j\left(\frac{\partial^2\eta^j}{\partial y^r
\,\partial y^s}+\frac{\partial\Gamma^j_{rs}}{\partial y^k}\,\eta^k+
\Gamma^j_{ks}\,\frac{\partial\eta^k}{\partial y^r}+\Gamma^j_{rk}\,
\frac{\partial\eta^k}{\partial y^s}\right).
\endgathered\hskip -2em
\tag2.3
\endgather
$$\par
     Note that there is an advanced theory of point symmetries for the
systems of differential equations (see \cite{7} or \cite{8}).
Application of this theory leads to the same determining equations
\thetag{2.2} and \thetag{2.3} for $\boldsymbol\eta$. Note also that 
\thetag{2.1} is a special form of point symmetry for equations \thetag{1.1}. 
Generic point symmetry has the form
$$
\boldsymbol\eta=\theta\,\frac{\partial}{\partial t}+
\xi\,\frac{\partial}{\partial x}+\sum^n_{i=1}\eta^i\,\frac{\partial}
{\partial y^i}.\hskip -2em
\tag2.4
$$
In this paper we restrict our consideration to the special point
symmetries \thetag{2.1} since they have transparent geometric
interpretation as vector fields on the manifold $M$. Geometric
interpretation of generic point symmetries \thetag{2.4} is the
subject for separate paper.\par
     Determining equations \thetag{2.2} and \thetag{2.3} for the field
of point symmetry \thetag{2.1} can be written without reference to
local coordinates. First of them provides vanishing of Lie derivative
of tensor field $\bold A$ (see more details in \cite{9}):
$$
L_{\boldsymbol\eta}({\bold A})=0.\hskip -2em
\tag2.5
$$
In case of $\det\bold A\neq 0$ from \thetag{2.2} and \thetag{2.3}
we derive the equation
$$
\sum^n_{k=1}\frac{\partial\eta^i}{\partial y^k}\,\Gamma^k_{rs}
=\frac{\partial^2\eta^i}{\partial y^r\,\partial y^s}+\sum^n_{k=1}
\left(\frac{\partial\Gamma^i_{rs}}{\partial y^k}\,\eta^k+
\Gamma^i_{ks}\,\frac{\partial\eta^k}{\partial y^r}+\Gamma^i_{rk}\,
\frac{\partial\eta^k}{\partial y^s}\right)\hskip -2em
\tag2.6
$$
that can replace \thetag{2.3}. If we take into account symmetry of
connection components $\Gamma^i_{rs}=\Gamma^i_{sr}$, which follows
from \thetag{1.1}, we can simplify this equation as below:
$$
\nabla_r\nabla_s\eta^i=\sum^n_{k=1}R^i_{srk}\,\eta^k.\hskip -2em
\tag2.7
$$
Here $R^i_{srk}$ are the components of curvature tensor, defined by
the formula
$$
R^i_{srk}=\frac{\partial\Gamma^i_{ks}}{\partial y^r}-\frac{\partial
\Gamma^i_{rs}}{\partial y^k}+\sum^n_{q=1}\Gamma^q_{ks}\,\Gamma^i_{rq}
-\sum^n_{q=1}\Gamma^q_{rs}\,\Gamma^i_{kq}\hskip -2em
\tag2.8
$$
(see for instance \cite{6}). In invariant (coordinateless) form the 
equations \thetag{2.7} are written as the relationship determining the 
commutator of Lie derivative $L_{\boldsymbol\eta}$ with covariant 
derivative along an arbitrary vector field $\bold Y$:
$$
[L_{\boldsymbol\eta},\,\nabla_{\bold Y}]=\nabla_{[\boldsymbol\eta,\,
\bold Y]}.\hskip -2em
\tag2.9
$$
Here one can remember the following lemma from \cite{10}.
\adjustfootnotemark{-1}
\proclaim{Lemma 2.1} Vector field $\boldsymbol\eta$ is an infinitesimal
affine transformation\footnote{here the term ``affine transformation''
is understood as a map preserving an affine connection (an automorphism
of an affine connection).} on the manifold $M$ if and only if the 
relationship
$$
[L_{\boldsymbol\eta},\,\nabla_{\bold Y}]\bold Z=\nabla_{[\boldsymbol\eta,\,
\bold Y]}\bold Z\hskip -2em
\tag2.10
$$
holds for two arbitrary vector fields $\bold Y$ and $\bold Z$ on this 
manifold.
\endproclaim
     The relationship $[L_{\boldsymbol\eta},\,\nabla_{\bold Y}]\varphi=
\nabla_{[\boldsymbol\eta,\,\bold Y]}\varphi$ for scalar field $\varphi$
holds identically. Therefore from \thetag{2.10} we obtain the
relationship 
$$
[L_{\boldsymbol\eta},\,\nabla_{\bold Y}]\bold W=\nabla_{[\boldsymbol\eta,\,
\bold Y]}\bold W\hskip -2em
\tag2.11
$$
which holds for arbitrary tensor field $\bold W$. It is equivalent to 
\thetag{2.9}. We state this result as a theorem.
\proclaim{Theorem 2.1} Vector field $\boldsymbol\eta$ is a field of point 
symmetry for the system of equations \thetag{1.1} with non-degenerate 
matrix $A$ if and only if it is an infinitesimal affine transformation for 
symmetric affine connection with components $\Gamma^j_{rs}$ and if Lie 
derivative of $\bold A$ along $\boldsymbol\eta$ is equal to zero (i\.~e\. 
both   relationships   \thetag{2.5}   and    \thetag{2.9}    hold 
simultaneously).
\endproclaim
\head
3. Operator field $\bold A$ and tensor fields relative to it.
\endhead
     According to the above geometric interpretation the system of 
equations \thetag{1.1} describes the dynamics in $x$ and $t$ for points 
of some manifold $M$ equipped with an affine connection $\Gamma$ and with 
operator field $\bold A$. Due to \thetag{2.5} operator field $\bold A$ is
invariant with respect to the field of point symmetry $\boldsymbol\eta$.
From \thetag{2.5} one can derive invariance for all integer powers of
the operator field $\bold A$. This means
$$
L_{\boldsymbol\eta}({\bold A^q})=0\text{, \ where \ }q\in\Bbb Z.
\hskip -2em
\tag3.1
$$
We suppose operator field $\bold A$ to be non-degenerate, therefore we 
admit negative integer powers $q<0$ in \thetag{3.1}.\par
     Besides with fields $\bold A^q$ the operator field $\bold A$ gives
rise to a set of tensor fields of type $(1,m)$, where $m\geqslant 1$ is 
a positive integer number. They are completely skew-symmetric in lower
indices. These fields are generated by integer powers of $\bold A$ and
Froelicher-Nijenhuis bracket from \cite{11}. Froelicher-Nijenhuis bracket 
plays an important role in the theory of equations of hydrodynamical type
$$
\frac{\partial y^i}{\partial\,t}=\sum^n_{j=1}A^i_j\,\frac{\partial y^j}
{\partial x},\quad i=1,\,\ldots,\,n,\hskip -2em
\tag3.2
$$
where manifold $M$ has no connection and the only differential geometric
object on it is the operator field $\bold A$. For instance in \cite{12}
we used Froelicher-Nijenhuis bracket for to obtain tensorial form $\bold P
=0$ of integrability condition for the equations \thetag{3.2} integrable
by means of generalized hodograph method from \cite{13}. Components of
tensor field $\bold P$ were expressed through components of matrix $A$.
Therefore $\bold P=0$ appears to be an effective integrability test,
which do not require the calculation of Riemann invariants required by
the theory from \cite{13}. Effective procedure for point classification
in the class of ordinary differential equations of the form
$$
y''=P(x,y)+3\,Q(x,y)\,y'+3\,R(x,y)\,{y'}^2+S(x,y)\,{y'}^3
$$
were derived in paper \cite{14}. However, this was done without use of 
Froelicher-Nijenhuis bracket.\par
     Let's consider tensor field $\bold B$ of type $(1,p)$ and tensor
field $\bold C$ of type $(1,q)$ made up by two vector fields $\bold b$ 
and $\bold c$ and two differential forms:
$$
\xalignat 2
&\bold B=\bold b\otimes\beta, &&\bold C=\bold c\otimes\gamma.
\hskip -2em
\tag3.3
\endxalignat
$$
For the fields \thetag{3.3} Froelicher-Nijenhuis bracket is calculated
as follows:
$$
\aligned
\{\bold B,\bold C\}&=[\bold b,\,\bold c]\otimes\beta\wedge
\gamma-\bold b\otimes L_{\bold c}\beta\wedge\gamma+
\bold c\otimes\beta\wedge L_{\bold b}\gamma+\\
\vspace{1ex}
&+(-1)^p\,\bold b\otimes \iota_{\bold c}\beta\wedge d\gamma +
(-1)^p\,\bold c\otimes d\beta\wedge\iota_{\bold b}\gamma.
\endaligned\hskip -2em
\tag3.4
$$
Here $d$ is an external differentiation defined by the formula
$$
\aligned
d\omega(\bold X_0,\ldots,\bold X_r)=\sum^r_{i=0}\frac{(-1)^i}{r+1}\,
\bold X_i(\omega(\bold X_0,\ldots,\hat\bold X_i,\ldots,\bold X_r))+\\
+\twosum_{0\leqslant i<j\leqslant r}\frac{(-1)^{i+j}}{r+1}\,
\omega([\bold X_i,\,\bold X_j],\bold X_0,\ldots,\hat\bold X_i,\ldots,
\hat\bold X_j,\ldots,\bold X_r),
\endaligned\hskip -2em
\tag3.5
$$
where $\omega$ is a differential $r$-form and $\bold X_0,\,\ldots,\,
\bold X_r$ are arbitrary vectorial fields. Hat over the sign of vector
field means that field with this particular number is omitted from the 
list of arguments of the form $\omega$. $L_{\bold b}$ and $L_{\bold c}$ 
in \thetag{3.4} are Lie derivatives along vector fields $\bold b$ and 
$\bold c$, while $\iota_{\bold b}$ and $\iota_{\bold c}$ are the 
differentiations of substitution. For $r$-form $\omega$ and vector field 
$\bold c$ the expression $\iota_{\bold c}\omega$ denotes $(r-1)$-form 
such that
$$
\iota_{\bold c}\omega(\bold X_1,\dots,\bold X_{r-1})=
r\,\omega(\bold c,\bold X_1,\dots,\bold X_{r-1}).\hskip -2em
\tag3.6
$$
The operation $\iota_{\bold c}$ is also known as {\it internal product} 
with respect to the vector field $\bold c$ (see \cite{9}). For the sake 
of completeness we give the formula that define external product $\beta
\wedge\gamma$ as a result of alternation of $\beta\otimes\gamma$:
$$
(\beta\wedge\gamma)(\bold X_1,\ldots,\bold X_{p+q})=
\sum_{\sigma}\frac{(-1)^\sigma}{(p+q)!}\,(\beta\otimes\gamma)(
\bold X_{\sigma_1},\ldots,\bold X_{\sigma(p+q)}).\hskip -2em
\tag3.7
$$
Here $\sigma$ is a permutation running over the whole group of 
permutations of the order $p+q$, while $(-1)^\sigma$ is a sign 
coefficient determined by the parity of $\sigma$.
\proclaim{Theorem 3.1} Froelicher-Nijenhuis bracket $\{\bold B,\bold C\}$ 
defined in \thetag{3.4} for tensor fields \thetag{3.3} can be expanded to 
the case of arbitrary two skew-symmetric tensor fields of types $(1,p)$ 
and $(1,q)$.
\endproclaim
     Proof of the theorem~3.1 based on formulas \thetag{3.4}, \thetag{3.5}, 
\thetag{3.6}, and \thetag{3.7} can be found in \cite{12}. On the base of
the same formulas one can check the following relationships:
$$
\aligned
&\{\bold B,\bold C\}+(-1)^{pq}\,\{\bold C,\bold B\}=0,\\
\vspace{1ex}
&(-1)^{rp}\,\{\{\bold B,\bold C\},\bold D\}+
(-1)^{pq}\,\{\{\bold C,\bold D\},\bold B\}+
(-1)^{qr}\,\{\{\bold D,\bold B\},\bold C\}=0.
\endaligned\hskip -3em
\tag3.8
$$
Here $\bold B$, $\bold C$, and $\bold D$ are completely skew-symmetric
tensor fields of types $(1,p)$, $(1,q)$, and $(1,r)$ respectively (such
fields are natural to call {\it vector valued differential forms}). Due
to the relationships \thetag{3.8} the set of vector valued differential
forms has a structure of graded Lie superalgebra over real numbers.
\proclaim{Theorem 3.2} Lie derivative $L_{\boldsymbol\eta}$ along 
arbitrary vector field $\boldsymbol\eta$ is a differentiation of degree
zero in Lie superalgebra of vector valued differential forms on $M$.
\endproclaim
    The statement of the theorem~3.2  is expressed by the following
equality being the Leibniz identity for Froelicher-Nijenhuis bracket:
$$
L_{\boldsymbol\eta}\{\bold B,\bold C\}=\{L_{\boldsymbol\eta}\bold B,
\bold C\}+\{\bold B,L_{\boldsymbol\eta}\bold C\}.\hskip -2em
\tag3.9
$$
For the fields \thetag{3.3} the identity \thetag{3.9} can be proved
by direct calculations based on the formula \thetag{3.4} and
commutational and anticommutational identities
$$
\xalignat 3
&L_{\boldsymbol\eta}\compos d-d\compos L_{\boldsymbol\eta}=0,
&&L_{\boldsymbol\eta}\compos\iota_{\boldsymbol\xi}-\iota_{\boldsymbol\xi}
\compos L_{\boldsymbol\eta}=\iota_{[\boldsymbol\eta,\,\boldsymbol\xi]},
&&\iota_{\boldsymbol\xi}\compos d+d\compos\iota_{\boldsymbol\xi}=
L_{\boldsymbol\xi}
\endxalignat
$$
from \cite{9}. For the case of arbitrary skew-symmetric tensor fields 
$\bold B$ and $\bold C$ of types $(1,p)$ and $(1,q)$ formula \thetag{3.9} 
is expanded by $\Bbb R$-linearity with the use of expansions
$$
\xalignat 2
&\bold B=\sum^{s_1}_{i=1}\bold b_i\otimes\beta_i,
&&\bold C=\sum^{s_2}_{j=1}\bold c_j\otimes\beta_j.
\endxalignat
$$\par\adjustfootnotemark{-1}
     Denote by $\Cal A=\langle\bold A^q,\ q\in\Bbb Z\,|\,\otimes,
C,\{*,*\}\rangle$ the linear span for the closure of the set of integer
powers of the operator field $\bold A$ with respect to the operations
of tensor product, contraction, and Froelicher-Nijenhuis 
bracket\footnote{permutation of indices are also assumed to be in the
list of operations generating algebra $\Cal A$. This means that $\Cal A$ 
is closed with respect to symmetrization and alternation of tensor 
fields.}. $\Cal A$ is $\Bbb R$-subalgebra in the algebra of tensor fields 
on $M$. For any two fields $\bold B$ and $\bold C$ from $\Cal A$ their 
Froelicher-Nijenhuis bracket $\{\bold B,\bold C\}$ (if it is defined) is 
an element of $\Cal A$. Tensor algebra $\Cal A$ generated by the field 
$\bold A$ contains Nijenhuis tensor $\bold N=\{\bold A,\bold A\}/2$ and 
Haantjes tensor from \cite{15} used in \cite{16} for to formulate 
diagonalizability criterion for the matrix $A$ in \thetag{3.2}. Moreover 
this algebra contains semihamiltonity tensor $\bold P$ (tensor of 
hydrodynamic integrability) constructed in \cite{12}. From \thetag{2.5} 
and \thetag{3.9} we derive the following theorem characterizing tensor 
algebra $\Cal A$.
\proclaim{Theorem 3.3} Suppose that $\bold B$ is an arbitrary tensor field
from the tensor algebra $\Cal A=\langle\bold A^q,\ q\in\Bbb Z\,|\,
\otimes,C,\{*,*\}\rangle$ generated by operator field $\bold A$ with 
components $A^i_j$ given by the system of equations \thetag{1.1}. Then
Lie derivative of $\bold B$ along point symmetry field $\boldsymbol\eta$ 
of the system \thetag{1.1} is zero, i\.~e\. $L_{\boldsymbol\eta}(\bold 
B)=0$.
\endproclaim
\head
4. Tensor fields generated by affine connection.
\endhead
    Affine connection $\Gamma$ is a second geometric structure on
$M$ given by the system of equations \thetag{1.1}. It defines tensor
field of curvature $\bold R$ with components \thetag{2.8}.
\proclaim{Theorem 4.1} Let $\bold R$ be the curvature tensor defined 
by affine connection $\Gamma$ from \thetag{1.1}. Then its Lie derivative 
along the field of point symmetry for the system of equations \thetag{1.1} 
is zero, i\.~e\. $L_{\boldsymbol\eta}(\bold B)=0$.
\endproclaim
\demo{Proof} In order to prove this theorem we shall use well known
commutational relationship for covariant derivatives (see \cite{17}):
$$
[\nabla_{\bold X},\,\nabla_{\bold Y}]\bold Z=
\nabla_{[\bold X,\,\bold Y]}\bold Z+\bold R(\bold X,\bold Y)\bold Z.\hskip 
-2em
\tag4.1
$$
Here $\bold X$, $\bold Y$, $\bold Z$ are three arbitrary vector fields on
$M$, while $\bold R(\bold X,\bold Y)$ is the operator field resulting by
contraction of curvature tensor $\bold R$ and vector fields $\bold X$ and 
$\bold Y$ with respect to last two indices in $\bold R$. Let's apply Lie
derivative $L_{\boldsymbol\eta}$ to both sides of the equality 
\thetag{4.1}. In left hand side we obtain
$$
L_{\boldsymbol\eta}[\nabla_{\bold X},\,\nabla_{\bold Y}]\bold Z=
L_{\boldsymbol\eta}\nabla_{\bold X}\nabla_{\bold Y}\bold Z-
L_{\boldsymbol\eta}\nabla_{\bold Y}\nabla_{\bold X}\bold Z.
\hskip -2em
\tag4.2
$$
Now in the expression $L_{\boldsymbol\eta}\nabla_{\bold X}
\nabla_{\bold Y}\bold Z$ we transpose Lie derivative $L_{\boldsymbol\eta}$
first with covariant derivative $\nabla_{\bold X}$, then with 
$\nabla_{\bold Y}$. As a result we have 
$$
\align
L_{\boldsymbol\eta}&\nabla_{\bold X}\nabla_{\bold Y}\bold Z=
[L_{\boldsymbol\eta},\,\nabla_{\bold X}]\nabla_{\bold Y}\bold Z
+\nabla_{\bold X}L_{\boldsymbol\eta}\nabla_{\bold Y}\bold Z=\\
\vspace{1ex}
&=[L_{\boldsymbol\eta},\,\nabla_{\bold X}]\nabla_{\bold Y}\bold Z
+\nabla_{\bold X}[L_{\boldsymbol\eta},\,\nabla_{\bold Y}]\bold Z
+\nabla_{\bold X}\nabla_{\bold Y}L_{\boldsymbol\eta}(\bold Z).
\endalign
$$
Taking into account \thetag{2.9} we come to the formula
$$
\pagebreak
L_{\boldsymbol\eta}\nabla_{\bold X}\nabla_{\bold Y}\bold Z
=\nabla_{[\boldsymbol\eta,\,\bold X]}\nabla_{\bold Y}\bold Z+
\nabla_{\bold X}\nabla_{[\boldsymbol\eta,\,\bold Y]}\bold Z
+\nabla_{\bold X}\nabla_{\bold Y}L_{\boldsymbol\eta}(\bold Z).
\hskip -2em
\tag4.3
$$
For $L_{\boldsymbol\eta}\nabla_{\bold Y}\nabla_{\bold X}\bold Z$
there is an analogous formula
$$
L_{\boldsymbol\eta}\nabla_{\bold Y}\nabla_{\bold X}\bold Z
=\nabla_{[\boldsymbol\eta,\,\bold Y]}\nabla_{\bold X}\bold Z+
\nabla_{\bold Y}\nabla_{[\boldsymbol\eta,\,\bold X]}\bold Z
+\nabla_{\bold Y}\nabla_{\bold X}L_{\boldsymbol\eta}(\bold Z).
\hskip -2em
\tag4.4
$$
Let's subtract \thetag{4.4} from \thetag{4.3} and let's apply the
identity \thetag{4.1}. As a result we can rewrite the equality
\thetag{4.2} as follows:
$$
\gathered
L_{\boldsymbol\eta}[\nabla_{\bold X},\,\nabla_{\bold Y}]\bold Z=
\nabla_{[[\boldsymbol\eta,\,\bold X],\,\bold Y]}\bold Z
+\bold R([\boldsymbol\eta,\,\bold X],\bold Y)\bold Z+
\nabla_{[\bold X,\,[\boldsymbol\eta,\,\bold Y]]}\bold Z\,+\\
\vspace{1ex}
+\,\bold R(\bold X,[\boldsymbol\eta,\,\bold Y])\bold Z+
\nabla_{[\bold X,\,\bold Y]}L_{\boldsymbol\eta}(\bold Z)+
\bold R(\bold X,\bold Y)L_{\boldsymbol\eta}(\bold Z).
\endgathered\hskip -2em
\tag4.5
$$
Applying Lie derivative $L_{\boldsymbol\eta}$ to last two summands
in right hand side of the identity \thetag{4.1} and taking into account
the commutational relationship \thetag{2.9} together with Jacobi
identity for commutator of vector fields we now obtain
$$
\align
&L_{\boldsymbol\eta}(\nabla_{[\bold X,\,\bold Y]}\bold Z)=
\nabla_{[[\boldsymbol\eta,\,\bold X],\,\bold Y]}\bold Z+
\nabla_{[\bold X,\,[\boldsymbol\eta,\,\bold Y]]}\bold Z+
\nabla_{[\bold X,\,\bold Y]}L_{\boldsymbol\eta}(\bold Z).
\hskip -2em
\tag4.6\\
\vspace{2ex}
&\aligned
L_{\boldsymbol\eta}(\bold R(\bold X,\bold Y)\bold Z)
=L_{\boldsymbol\eta}(\bold R&)(\bold X,\bold Y)\bold Z
+\bold R([\boldsymbol\eta,\,\bold X],\bold Y)\bold Z\,+\\
\vspace{1ex}
&+\,\bold R(\bold X,[\boldsymbol\eta,\,\bold Y])\bold Z
+\bold R(\bold X,\bold Y)L_{\boldsymbol\eta}(\bold Z).
\endaligned\hskip -2em
\tag4.7
\endalign
$$
When subtracting \thetag{4.6} and \thetag{4.7} from \thetag{4.5} all
summands in right hand side are canceled, except for only one. Here
is the ultimate result of the above calculations:
$$
L_{\boldsymbol\eta}(\bold R)(\bold X,\bold Y)\bold Z=0.
$$
Since $\bold X$, $\bold Y$, and $\bold Z$ are arbitrary vector fields,
we obtain the equality $L_{\boldsymbol\eta}(\bold R)=0$, which was
required to prove the theorem.\qed\enddemo
     Covariant derivative $\nabla_{\bold Y}$ is closely connected with
the concept of {\it covariant differential}. For the tensor field 
$\bold W$ of type $(r,s)$ its covariant differential is a tensor field
$\nabla\bold W$ of type $(r,s+1)$ such that when being contracted with
vector field $\bold Y$ yields the field of covariant derivative 
$\nabla_{\bold Y}\bold W$:
$$
\nabla_{\bold Y}\bold W=C(\bold Y\otimes\nabla\bold W).\hskip -2em
\tag4.8
$$
Let's apply Lie derivative $L_{\boldsymbol\eta}$ along the field of 
point symmetry of \thetag{1.1} to both sides of the equality \thetag{4.8}.
And take into account \thetag{2.9} doing this:
$$
\align
&\aligned
L_{\boldsymbol\eta}(\nabla_{\bold Y}\bold W)=&\nabla_{[\boldsymbol\eta,
\,\bold Y]}\bold W+\nabla_{\bold Y}L_{\boldsymbol\eta}(\bold W)=\\
\vspace{1ex}
&\ =C([\boldsymbol\eta,\,\bold Y]\otimes\nabla\bold W)+C(\bold Y\otimes
\nabla L_{\boldsymbol\eta}(\bold W)).
\endaligned\\
\vspace{2ex}
&L_{\boldsymbol\eta}C(\bold Y\otimes\nabla\bold W)=
C([\boldsymbol\eta,\,\bold Y]\otimes\nabla\bold W)+
C(\bold Y\otimes L_{\boldsymbol\eta}(\nabla\bold W)).
\endalign
$$
Equating right hand sides of two resulting relationships we obtain
$$
\pagebreak 
C(\bold Y\otimes\nabla L_{\boldsymbol\eta}(\bold W))=
C(\bold Y\otimes L_{\boldsymbol\eta}(\nabla\bold W)).
$$
Now remember that $\bold Y$ is an arbitrary vector field. This leads
us to the equality $\nabla L_{\boldsymbol\eta}(\bold W)=L_{\boldsymbol\eta}
(\nabla\bold W)$, which can be written as
$$
[L_{\boldsymbol\eta},\,\nabla]=0\hskip -2em
\tag4.9
$$
since $\bold W$ is also an arbitrary tensor field.
\proclaim{Theorem 4.2} Lie derivative $L_{\boldsymbol\eta}$ along the field 
of point symmetry for the system of equations \thetag{1.1} is commutating 
with covariant differential $\nabla$ defined by the affine connection with
components $\Gamma^j_{rs}$ from \thetag{1.1}.
\endproclaim
\adjustfootnotemark{-1}
     Let's add the curvature tensor to the set of integer powers of the
operator field $\bold A$ and let's consider linear span for the closure
of the resulting set with respect to the operations of tensor product, 
contraction, Froelicher-Nijenhuis bracket, and covariant 
differential\footnote{permutation of indices here are assumed to be in
the set of generating operations for the algebra $\Cal R$ as well as in 
case of algebra $\Cal A$.}:
$$
\Cal R=\langle\bold R,\,\bold A^q,\ q\in\Bbb Z\,|\,\otimes,C,\{*,*\},
\nabla\rangle.\hskip -2em
\tag4.10
$$
\proclaim{Theorem 4.3} If tensor algebra $\Cal R$ in \thetag{4.10} is
generated by the operator field $\bold A$ and affine connection $\Gamma$
with components taken from \thetag{1.1}, then for any tensor field 
$\bold B$ from $\Cal R$ its Lie derivative $L_{\boldsymbol\eta}(\bold B)$
along the field of point symmetry $\boldsymbol\eta$ for this system of 
equations is zero, i\.~e\. $L_{\boldsymbol\eta}(\bold B)=0$.
\endproclaim
     Tensor algebra $\Cal R$ contains two tensor fields important for the
further study of the equations \thetag{1.1}. These are the fields $\bold R$ 
and $\bold S$ with components
$$
\xalignat 2
&R_{ij}=\sum^n_{k=1}R^k_{ikj},&&S_{ij}=\sum^n_{k=1}R^k_{kij}.\hskip -2em
\tag4.11
\endxalignat
$$

First of them is a field of Ricci tensor. For arbitrary affine connection
it has symmetric part and antisymmetric part as well:
$$
\xalignat 2
&\hat R_{ij}=\frac{R_{ij}+R_{ji}}{2},
&&\tilde R_{ij}=\frac{R_{ij}-R_{ji}}{2}.
\hskip -2em
\tag4.12
\endxalignat
$$
Tensor fields $\hat\bold R$ and $\tilde\bold R$ with components 
\thetag{4.12} belong to the tensor algebra $\Cal R$. For $\hat\bold R$ 
and $\tilde\bold R$ we have the relation $\bold S=2\,\tilde\bold R$, 
which follows from well-known identity $R^s_{ijk}+R^s_{jki}+R^s_{kij}
=0$. The latter is the consequence of symmetry $\Gamma^j_{rs}=
\Gamma^j_{sr}$.\par
\head
5. Symmetries with stationary point.
\endhead
     Let $\boldsymbol\eta$ be a field of point symmetry for the system
of equations \thetag{1.1}. It gives rise to the local one-parametric group
of transformations $\varphi_\tau\!:M\to M$. In local coordinates $y^1,\,
\ldots,\,y^n$ on $M$ these transformations are represented by functions
$$
\align
x^1&=x^1(\tau,y^1,\ldots,y^n),\hskip -2em\\
.\ .&\ .\ .\ .\ .\ .\ .\ .\ .\ .\ .\ .\ .\ .\ .\ \hskip -2em
\tag5.1\\
x^n&=x^n(\tau,y^1,\ldots,y^n).\hskip -2em
\endalign
$$
Here arguments $y^1,\,\ldots,\,y^n$ are the coordinates of initial point
$p\in M$, while the values of functions $x^1,\,\ldots,\,x^n$ are the
coordinates of transformed point $\varphi_\tau(p)$. Functions \thetag{5.1} 
are the solutions for the system of ordinary differential equations
$$
\align
(x^1)'_\tau&=\eta^1(x^1,\ldots,x^n),\hskip -2em\\
.\ .\ .\ .&\ .\ .\ .\ .\ .\ .\ .\ .\ .\ .\ .\ .\ .\ \hskip -2em
\tag5.2\\
(x^n)'_\tau&=\eta^n(x^1,\ldots,x^n)\hskip -2em
\endalign
$$
such that they solve Cauchy problem with the following initial data:
$$
x^1(\tau)\,\hbox{\vrule height 8pt depth 8pt width 0.5pt}_{\,\tau=0}=y^1,
\ \ldots,\ x^n(\tau)\,\hbox{\vrule height 8pt depth 8pt width 0.5pt}_{\,
\tau=0}=y^n.\hskip -2em
\tag5.3
$$
Initial data \thetag{5.3} mean that $\varphi_\tau\!:M\to M$ is identical
map for $\tau=0$.\par
     Let $p_0$ be a stationary point for the maps $\varphi_\tau$, i\.~e\.
$\varphi_\tau(p_0)=p_0$ for all $\tau$. If $y^1_0,\,\ldots,\,y^n_0$ are
coordinates of such point, then functions
$$
x^1(\tau)=y^1_0=\const,\ \ldots,\ x^n(\tau)=y^n_0=\const
$$
should satisfy the system of equations \thetag{5.2} and initial conditions
\thetag{5.3}. Substitution of these functions into \thetag{5.2} leads us
to the following well-known fact.
\proclaim{Lemma 5.1} Point $p$ on $M$ is a stationary point of local 
one-parametric group of transformations $\varphi_\tau\!:M\to M$ if and
only if corresponding vector field $\boldsymbol\eta$ vanishes at this
point.
\endproclaim
     If $p_0$ is a vanishing point for the vector field $\boldsymbol\eta$,
then we can consider the differential of the map $\varphi_\tau$. Due to
$\varphi_\tau(p_0)=p_0$ it is linear operator in tangent space $T_{p_0}(M)$ 
to the manifold $M$ at the point $p_0$. In local coordinates $y^1,\,\ldots,
\,y^n$ operator $\varphi_{\tau*}\!:T_{p_0}(M)\to T_{p_0}(M)$ is expressed 
by the Jacobi matrix for the system of functions \thetag{5.1} at 
stationary point $p_0$:
$$
\Phi(\tau)=\Vmatrix \dfrac{\partial x^1}{\partial y^1} & \hdots &
\dfrac{\partial x^1}{\partial y^n}\\ \vspace{-2pt}\vdots & \ddots
& \vdots\\ \vspace{4pt}\dfrac{\partial x^n}{\partial y^1} & \hdots
& \dfrac{\partial x^n}{\partial y^n}\endVmatrix\hskip -2em
\tag5.4
$$
Components of matrix \thetag{5.4} are smooth functions of parameter
$\tau$. Denote by $F$ the matrix composed by derivatives in $\tau$ of 
the components of matrix $\Phi(\tau)$ for $\tau=0$. In order to find
the matrix $\Phi(\tau)$ we should solve Cauchy problem \thetag{5.3} for 
the equations \thetag{5.2}. However, components of the matrix $F$
can be expressed through components of vector field $\boldsymbol\eta$
without solving the equations \thetag{5.2}:
$$
\pagebreak
F^i_j=\frac{\partial\Phi^i_j}{\partial\tau}\,\hbox{\vrule height 14pt
depth 8pt width 0.5pt}_{\,\tau=0}=\frac{\partial\eta^i}{\partial y^j}
\,\hbox{\vrule height 14pt depth 8pt width 0.5pt}_{\,p=p_0}.\hskip -2em
\tag5.5
$$\par
     Matrices $\Phi(\tau)$ with various values of $\tau$ form one-parametric
matrix group because $\Phi(\tau_1+\tau_2)=\Phi(\tau_1)\,\Phi(\tau_2)$. These
matrices are analogs of rotation matrices, while matrix $F$ with components
\thetag{5.5} is an analog of infinitesimal rotation around the point $p_0$. 
From $\Phi(\tau_1+\tau_2)=\Phi(\tau_1)\,\Phi(\tau_2)$ we derive differential 
equation for $\Phi(\tau)$:
$$
\Phi'_\tau=F\,\Phi.\hskip -2em
\tag5.6
$$
The solution of \thetag{5.6} is a matrix exponent $\Phi(\tau)=\exp(F\,
\tau)$.\par
     We can give another (more direct) interpretation for the matrix $F$.
Let's consider Lie derivative $L_{\boldsymbol\eta}$. Being a 
differentiation in the set of vector fields Lie derivative $L_{\boldsymbol
\eta}$ satisfies the following relationships:
$$
\xalignat 2
&L_{\boldsymbol\eta}(\bold X+\bold Y)=L_{\boldsymbol\eta}(\bold X)+
L_{\boldsymbol\eta}(\bold Y),&&L_{\boldsymbol\eta}(\psi\cdot\bold X)=
\psi\cdot L_{\boldsymbol\eta}(\bold X)+(\boldsymbol\eta\psi)\cdot\bold X.
\endxalignat
$$
Here $\bold X$ and $\bold Y$ are vector fields, while $\psi$ is a scalar
field. When written at the point $p_0$, due to $\boldsymbol\eta(p_0)=0$,
the above relationships look like $L_{\boldsymbol\eta}(\bold X+\bold Y)=
L_{\boldsymbol\eta}(\bold X)+L_{\boldsymbol\eta}(\bold Y)$ and 
$L_{\boldsymbol\eta}(\psi\cdot\bold X)=\psi\cdot L_{\boldsymbol\eta}
(\bold X)$. 
This means that $L_{\boldsymbol\eta}$ acts as a linear operator
$$
L_{\boldsymbol\eta}\!:T_{p_0}(M)\to T_{p_0}(M)\hskip -2em
\tag5.7
$$
in the tangent space to $M$ at the point $p_0$. By direct calculations
in local coordinates one can easily check that the matrix of the 
operator \thetag{5.7} has the components defined by the right hand side 
of \thetag{5.5}.
\head
6. Estimate for dimension of symmetry algebra. 
\endhead
     Vector fields of point symmetries for the system of equations
\thetag{1.1} constitute Lie algebra over real numbers (see for instance
\cite{7} or \cite{8}). In this paper we consider a subalgebra that
consists of special point symmetries \thetag{2.1}, denote it by
$\Cal G$. Elements of $\Cal G$ are interpreted as vector fields on the
manifold $M$. Let $\boldsymbol\eta\in\Cal G$. Consider the equation
\thetag{4.9}, which is equivalent to \thetag{2.9}. In local coordinates
$y^1,\,\ldots,\,y^n$ it is written in form of the system of equations
\thetag{2.6} with respect to the components of vector fields $\boldsymbol
\eta$. Let's rewrite this system of equations as 
$$
\frac{\partial F^i_s}{\partial y^r}=\sum^n_{k=1}\Gamma^k_{rs}\,F^i_k
-\sum^n_{k=1}\left(\frac{\partial\Gamma^i_{rs}}{\partial y^k}\,\eta^k+
\Gamma^i_{ks}\,F^k_r+\Gamma^i_{rk}\,F^k_s\right)\hskip -2em
\tag6.1
$$
and let's complete it with the equations defining $F^i_r$:
$$
\frac{\partial\eta^i}{\partial y^r}=F^i_r.\hskip -2em
\tag6.2
$$
Equations \thetag{6.1} and \thetag{6.2} in the aggregate form complete
system of Pfaff equations with respect to $n^2+n$ functions $F^i_s(y^1,
\ldots,y^n)$ and $\eta^i(y^1,\ldots,y^n)$. It's known that each solution
of such system is uniquely defined by the initial data at some point:
$$
\pagebreak 
\xalignat 2
&F^i_s\,\hbox{\vrule height 8pt depth 8pt width 0.5pt}_{\,p=p_0}=
F^i_s(0),&&\eta^i\,\hbox{\vrule height 8pt depth 8pt width
0.5pt}_{\,p=p_0}=\eta^i(0).\hskip -2em
\tag6.3
\endxalignat
$$
However, not for each initial data \thetag{6.3} one can find appropriate 
solution of the equations \thetag{6.1} and \thetag{6.2}. Some constraints
there appear since system o Pfaff equations my be not completely 
compatible. Therefore we have the theorem.
\proclaim{Theorem 6.1} For any system of equations \thetag{1.1} the Lie
algebra of special point symmetries \thetag{2.1} is finite dimensional
and its dimension is not greater than $n(n+1)$.
\endproclaim
     By proving the above theorem~6.1 we reproduce partially (on a 
level of Lie algebras) the following well-known result (see \cite{10}, 
chapter~4, \S\,1).
\proclaim{Theorem 6.2} Let $M$ be $n$-dimensional manifold with affine 
connection. Then group of affine transformations of the manifold $M$ is 
a Lie group of the dimension not greater than $n(n+1)$.
\endproclaim
\head
7. Case of maximal degeneration.
\endhead
     The estimate in theorem~6.1 is exact. For to prove this we shall
find the condition for complete compatibility of the system of Pfaff
equations \thetag{6.1} and \thetag{6.2}. For the sake of brevity let's
introduce the notations
$$
\Omega^{ij}_{rsq}=\Gamma^j_{rs}\,\delta^i_q-\Gamma^i_{qs}\,\delta^j_r
-\Gamma^i_{qr}\,\delta^j_s,\hskip -2em
\tag7.1
$$
where $\delta^i_j$ is a Kronecker's delta-symbol. Then we can write
\thetag{6.1} and \thetag{6.2} as 
$$
\xalignat 2
&\quad\frac{\partial F^i_s}{\partial y^r}=\sum^n_{j=1}\sum^n_{q=1}
\Omega^{ij}_{rsq}\,F^q_j-\sum^n_{k=1}\frac{\partial\Gamma^i_{rs}}
{\partial y^k}\,\eta^k,&&\frac{\partial\eta^i}{\partial y^r}=F^i_r.
\hskip -2em
\tag7.2
\endxalignat
$$
Compatibility condition for the equations \thetag{7.2} is obtained
by equating partial derivatives calculated by virtue of these
equations:
$$
\xalignat 2
&\quad\frac{\partial^2 F^i_s}{\partial y^p\,\partial y^r}=\frac{\partial^2
F^i_s}{\partial y^r\,\partial y^p},&&\frac{\partial^2\eta^i}{\partial
y^p\,\partial y^r}=\frac{\partial^2\eta^i}{\partial y^r\,\partial y^p}.
\hskip -2em
\tag7.3
\endxalignat
$$
Note that second part of equalities \thetag{7.3} are identically fulfilled
due to the symmetry $\Gamma^i_{pr}=\Gamma^i_{rp}$ and $\Omega^{ij}_{rpq}=
\Omega^{ij}_{prq}$. For to write complete compatibility condition we
should equate coefficients of each $F^q_j$ and each $\eta^k$ in the first
part of \thetag{7.3} assuming $F^q_j$ and $\eta^k$ to be independent 
variables. This yields
$$
\gather
\aligned
\frac{\partial\Omega^{ij}_{rsq}}{\partial y^p}+&\sum^n_{\alpha=1}
\sum^n_{\beta=1}\Omega^{i\beta}_{rs\alpha}\,\Omega^{\alpha j}_{p\beta q}
-\frac{\partial\Gamma^i_{rs}}{\partial y^q}\,\delta^j_p=\\
&=\frac{\partial\Omega^{ij}_{psq}}{\partial y^r}+\sum^n_{\alpha=1}
\sum^n_{\beta=1}\Omega^{i\beta}_{ps\alpha}\,\Omega^{\alpha j}_{r\beta q}
-\frac{\partial\Gamma^i_{ps}}{\partial y^q}\,\delta^j_r,
\endaligned\hskip -2em
\tag7.4\\
\vspace{2ex}
\frac{\partial^2\Gamma^i_{rs}}{\partial y^p\,\partial y^k}+\sum^n_{q=1}
\sum^n_{j=1}\Omega^{ij}_{rsq}\,\frac{\partial\Gamma^q_{pj}}{\partial y^k}=
\frac{\partial^2\Gamma^i_{ps}}{\partial y^r\,\partial y^k}+\sum^n_{q=1}
\sum^n_{j=1}\Omega^{ij}_{psq}\,\frac{\partial\Gamma^q_{rj}}{\partial y^k}.
\hskip -2em
\tag7.5
\endgather
$$
Let's substitute \thetag{7.1} into \thetag{7.4} and do contract indices
$j$ and $q$. As a result of this operation we obtain the equality
$$
\frac{\partial\Gamma^i_{rs}}{\partial y^p}+\sum^n_{\beta=1}
\Gamma^\beta_{rs}\,\Gamma^i_{p\beta}=\frac{\partial\Gamma^i_{ps}}
{\partial y^r}+\sum^n_{\beta=1}\Gamma^\beta_{ps}\,\Gamma^i_{r\beta},
\hskip -2em
\tag7.6
$$
which means $\bold R=0$ (compare with \thetag{2.8}). Upon substituting
\thetag{7.1} into \thetag{7.5} this relationship appears to be
differential consequence of \thetag{7.6}: it can be derived from
\thetag{7.6} by differentiating both sides with respect to $y^k$.
As for the relationship \thetag{7.4} itself, after substituting
\thetag{7.1} into it and after collecting similar term it takes
the form of linear combination of four equalities, each obtained
by some permutation of indices in \thetag{7.6}. The result
of these computations we formulate as a lemma.
\proclaim{Lemma 7.1} System of Pfaff equations \thetag{6.1} and 
\thetag{6.2} is completely compatible if and only if affine
connection with components $\Gamma^i_{rs}$ is a flat connection
with zero curvature tensor.
\endproclaim
    The following fact is well-known in geometry (see \cite{17}):
if affine connection is flat, then there exist a local coordinate 
system $y^1,\,\ldots,\,y^n$ on $M$ such that all connection 
components are identically zero in these coordinates (they are called 
{\it euclidean coordinates}). We shall use such euclidean coordinates 
in order to solve Pfaff equations \thetag{6.1} and \thetag{6.2}. 
Here these equations become very simple:
$$
\xalignat 2
&\frac{\partial F^i_s}{\partial y^r}=0,
&&\frac{\partial\eta^i}{\partial y^r}=F^i_r.\hskip -2em
\tag7.7
\endxalignat
$$
Cauchy problem \thetag{6.3} for the equations \thetag{7.7} is solvable
for any initial data $F^i_s(0)$ and $\eta^i(0)$. If we suppose the
coordinates of $p_0$ to be zero, then the solution of this Cauchy
problem is given by the following linear functions:
$$
\xalignat 2
&F^i_s=F^i_s(0)=\const,&&\eta^i=\eta^i(0)+\sum^n_{s=1}F^i_s(0)
\cdot y^s.\hskip -2em
\tag7.8
\endxalignat
$$\par
    Taking components of the field $\boldsymbol\eta$ from \thetag{7.8}
we substitute this field into the equation \thetag{2.5}. When 
written in local coordinates, this equation has the form \thetag{2.2}. 
As a result we obtain the equation
$$
\sum^n_{k=1}\left(A^i_k\,F^k_j(0)-F^i_k(0)\,A^k_j\right)+
\sum^n_{k=1}\!\left(\,\shave{\sum^n_{s=1}}F^k_s(0)\,y^s+\eta^k(0)
\!\right)\frac{\partial A^i_j}{\partial y^k}=0.\hskip -3em
\tag7.9
$$
For the arbitrary vector field $\boldsymbol\eta$ with components 
\thetag{7.8} to be the field of point symmetry for \thetag{1.1}
the equations \thetag{7.9} should be fulfilled identically for
arbitrary values of constants $F^i_s(0)$ and $\eta^i(0)$. Let's
substitute $F^i_s(0)=0$ into \thetag{7.9} and take into account
that $\eta^i(0)$ are arbitrary constants. This yields $A^i_j=\const$
and
$$
\pagebreak
\sum^n_{k=1}A^i_k\,F^k_j(0)=\sum^n_{k=1}F^i_k(0)\,A^k_j.\hskip -3em
\tag7.10
$$
The relationships \thetag{7.10} mean that matrix $A$ should commutate
with arbitrary matrix $F$ whose components are $F^i_j(0)$. Therefore
$A$ can differ from unit matrix only by scalar factor.
\proclaim{Theorem 7.1} System of equation \thetag{1.1} with nondegenerate 
matrix $A$ has an algebra of point symmetries \thetag{2.1} of maximal
dimension $n(n+1)$ if and only if $A$ is a scalar matrix ($A^i_j=a\,
\delta^i_j$) with constant factor $a=\const$ and if affine connection
$\Gamma$ with components $\Gamma^j_{rs}$ is flat (i\.~e\. its curvature
is zero).
\endproclaim
     This result is in agreement with the following well-known 
geometric fact (see \cite{10}, chapter~4, \S\,1).
\proclaim{Theorem 7.2} Let $M$ be $n$-dimensional manifold with
affine connection. The dimension of the group of affine transformations
is equal to $n(n+1)$ exactly if and only if $M$ is an ordinary flat
affine space $\Bbb A_n$ with natural flat affine connection.
\endproclaim
     Systems of equations \thetag{1.1} described by the theorem~7.1
constitute maximally symmetric, and therefore maximally simple (maximally 
degenerate), subclass of such systems. For any of them we can find point
transformation \thetag{1.4} breaking this system into separate equations
of the form 
$$
\frac{\partial y^i}{\partial\,t}=a\,\frac{\partial^2 y^i}{\partial\,x^2},
\quad i=1,\,\ldots,\,n,\hskip -2em
\tag7.11
$$
where $a$ is a common constant. Geometry associated with the equations
\thetag{7.11} is also maximally simple. In order to find the equations
with more interesting geometry further we consider the cases when
algebra of point symmetries of such equations is not maximal, but is
extremal in some sense as described below.
\head
8. Cases of intermediate degeneration.
\endhead
     Let $p_0$ be some fixed point on the manifold $M$. Vector fields
from $\Cal G$ vanishing at the point $p_0$ constitute subalgebra in
$\Cal G$. Denote it by $\Cal G(p_0)$. For the dimension of factorspace 
$\Cal G/\Cal G(p_0)$ (not factoralgebra) we have an obvious estimate 
$$
\dim(\Cal G/\Cal G(p_0))\leqslant\dim M=n.\hskip -2em
\tag8.1
$$
Estimate \thetag{8.1} can be obtained from the statement of Cauchy 
problem \thetag{6.3} for Pfaff equations \thetag{6.1} and \thetag{6.2}.
\par
     When Pfaff equations \thetag{6.1} and \thetag{6.2} are not completely
compatible the corresponding Cauchy problem for them can be solved not for
all initial data in \thetag{6.3}. Restrictions are due to differential
consequences of \thetag{6.1} and \thetag{6.2}. One of such differential
consequences can be formulated in terms of symmetric part of Ricci tensor: 
$L_{\boldsymbol\eta}(\hat\bold R)=0$. Let's write this equation
in local coordinates for the field of point symmetry $\boldsymbol\eta\in
\Cal G(p_0)$:
$$
\pagebreak
\sum^n_{k=1}\hat R_{ik}\,F^k_j(0)+\sum^n_{k=1}\hat R_{kj}\,F^k_i(0)=0.
\hskip -2em
\tag8.2
$$
Here $F^i_j(0)$ are the components of the matrix defined in \thetag{5.5}.
These quantities are initial data in Cauchy problem \thetag{6.3}.\par
     Ricci tensor $\hat\bold R$ is symmetric tensor field of type $(0,2)$,
i\.~e\. it is a field of quadratic forms. Denote by $\hat R(\bold X,\bold 
Y)$ the symmetric bilinear form given by the field $\hat\bold R$. Let $m$ 
be the sum of positive and negative indices in signature of this form
(see for instance \cite{18}). In other words $m$ is the rank of matrix
formed by the components of Ricci tensor $\hat\bold R$.
\definition{Definition 8.1} Say that the system of equations \thetag{1.1}
belongs to $m$-th case of intermediate degeneration ($1\leqslant m\leqslant 
n-1$) if the rank of symmetric part of Ricci tensor $\hat\bold R$ defined
by affine connection $\Gamma$ with components $\Gamma^j_{rs}$ is equal to
$m$ everywhere on the manifold $M$.
\enddefinition
    Case of maximal degeneration corresponds to $m=0$, while in case of
general position we have $m=n$. In order to find an estimate for the 
dimension of algebra of point symmetries depending on $m$ we write 
\thetag{8.2} in invariant form
$$
\hat R(\bold F\bold X,\bold Y)=-\hat R(\bold X,\bold F\bold Y),
\hskip -2em
\tag8.3
$$
which do not refer to local coordinates. This equality should be fulfilled 
for two arbitrary vectors $\bold X$ and $\bold Y$ from tangent space 
$V=T_{p_0}(M)$. Here $\bold F=L_{\boldsymbol\eta}$ is a linear operator
from  \thetag{5.7}. Denote by $W$ kernel of quadratic form $\hat R$ (see 
\cite{18}):
$$
W=\{\bold X\in V:\ \ \hat R(\bold X,\bold Y)=0\text{\ \ for all \ }
\bold Y\in V\}.\hskip -2em
\tag8.4
$$
Its easy to calculate the dimension of the subspace \thetag{8.4}: $\dim W
=n-m$.\par 
     From \thetag{8.3} we find that $W$ is invariant under the action
of the operator $\bold F$, i\.~e\. $\bold X\in W$ implies $\bold F\bold X
\in W$. Operators preserving $W$ constitute a subspace of the dimension 
$n^2-m(n-m)$ in the space $\End(V)$ of all linear operators in $V$. Denote 
it by $\End(V|W)$. For the operators from $\End(V|W)$ the concept of {\it 
restriction to $W$} and the concept of {\it factoroperator} are defined:
$$
\xalignat 2
&\bold F\,\hbox{\vrule height 8pt depth 8pt width 0.5pt}_{\,W}: W\to W,
&&\hat\bold F=\bold F\,\hbox{\vrule height 8pt depth 8pt width 
0.5pt}_{\,V/W}:V/W\to V/W.\hskip -2em
\tag8.5
\endxalignat
$$
The action of latter one to cosets relative to $W$ is given by the formula 
$$
\hat\bold F\Cl_{\sssize W}(\bold X)=\Cl_{\sssize W}(\bold F\bold X)
$$ 
(see more details in \cite{18}). Replacing operator $\bold F\in\End(V|W)$ 
by factoroperator $\hat\bold F\in\End(V/W)$ we factorize over the operators 
mapping $V$ into $W$:
$$
\End(V/W)\cong\End(V|W)/\Hom(V,W).\hskip -2em
\tag8.6
$$
Subspace $W$ is a kernel of $\hat R$. This implies that $\hat R$ induces
nondegenerate symmetric bilinear form $\hat R$ in factorspace $V/W$. It 
is defined by formula
$$
\hat R(\hat\bold X,\hat\bold Y)=\hat R(\bold X,\bold Y),\hskip -2em
\tag8.7
$$
where $\hat\bold X=\Cl_{\sssize W}(\bold X)$ and $\hat\bold Y=\Cl_{
\sssize W}(\bold Y)$ are cosets relative to the subspace $W$. In whole,
the equality \thetag{8.3} for the operators $\bold F\in\End(V)$ appears
to be equivalent to the following two conditions:
\roster
\item operator $\bold F$ belongs to the subspace $\End(V|W)$ that consists
      of operators having $W=\Ker\hat R$ as invariant subspace;
\item factoroperator $\hat F$ in \thetag{8.5} is skew-symmetric with
      respect to the symmetric bilinear form \thetag{8.7}, i\.~e\. the
      equality $\hat R(\hat\bold F\hat\bold X,\hat\bold Y)=-\hat R(\hat
      \bold X,\hat\bold F\hat\bold Y)$ holds for two arbitrary vectors
      $\hat\bold X$ and $\hat\bold Y$ from factorspace $V/W$.
\endroster
The dimension of subspace formed by the operators skew-symmetric with
respect to nondegenerate bilinear form \thetag{8.7} is given by the
formula
$$
\dim\hat\Cal F_{\sssize\text{skew}}=\frac{m(m-1)}{2}.
$$
We know the dimension of the space $\Hom(V,W)$, over which we factorize
in \thetag{8.6}. Thereby we can calculate the dimension of the space of
operators satisfying \thetag{8.3}:
$$
\dim\Cal F=\dim\Hom(V,W)+\dim\hat\Cal F_{\sssize\text{skew}}=
n(n-m)+\frac{m(m-1)}{2}.\hskip -3em
\tag8.8
$$
Let's combine \thetag{8.8} with \thetag{8.1}. This gives an estimate 
for dimensions of symmetry algebras of the system of equations 
\thetag{1.1} in all cases of intermediate degeneration:
$$
\dim(\Cal G)\leqslant n(n+1-m)+\frac{m(m-1)}{2}.\hskip -2em
\tag8.9
$$
This estimate is in agreement with the following theorem from \cite{19}.
\proclaim{Theorem 8.1} Maximally movable spaces of affine connection
with the symmetric part of Ricci tensor of the rank $m$ possess 
transitive groups of automorphisms with $n(n+1-m)+m(m-1)/2$ parameters.
\endproclaim
    In paper \cite{20} was shown that the estimate \thetag{8.9} is exact, 
some examples of spaces, where the upper bound is reached, were given. Here 
we give complete description of such spaces and describe the appropriate 
equations \thetag{1.1} for them in each case of intermediate degeneration. 
For $m=1$ this was done in \cite{21}.
\head
9. Structure of the tensor field $\bold S$.
\endhead
     Let \thetag{1.1} be the system of equations belonging to $m$-th 
case of intermediate degeneration such that its algebra of point 
symmetries has maximal dimension
$$
\dim(\Cal G)=n(n+1-m)+\frac{m(m-1)}{2}.\hskip -2em
\tag9.1
$$
Let $\bold S=2\,\tilde\bold R$ be a tensor field of type $(0,2)$ with
components \thetag{4.11}. For this field we state the following lemma.
\proclaim{Lemma 9.1} Kernel of skew-symmetric bilinear form defined by
$\bold S$ contains the kernel of the form defined by symmetric part of 
Ricci tensor.
\endproclaim
\demo{Proof} Tensor field $\bold S$ belongs to the algebra $\Cal R$.
Therefore due to theorem~4.3 we have the relationship $L_{\boldsymbol
\eta}\bold S=0$ satisfied for any field of point symmetry $\boldsymbol
\eta$ from $\Cal G$. Repeating arguments from section~8 we get
the relationship
$$
S(\bold F\bold X,\bold Y)=-S(\bold X,\bold F\bold Y)\hskip -2em
\tag9.2 
$$
analogous to \thetag{8.3}. The equality \thetag{9.1} here means that
\thetag{9.2} should be fulfilled for all operators satisfying \thetag{8.3}. 
Let's take an operator $\bold F\in\Hom(V,W)$ of special form 
$\bold F=\bold w\otimes\alpha$, where $\bold w$ is a vector from the kernel
of the form $\hat R$ and $\alpha$ is some arbitrary linear functional in
the space $V$. Operator $\bold F$ satisfies the equation \thetag{8.3}.
Therefore it should satisfy the equation \thetag{9.2} too. Substituting
$\bold F=\bold w\otimes\alpha$ into \thetag{9.2} we get the following 
equality:
$$
\alpha(\bold X)\,S(\bold w,\bold Y)=\alpha(\bold Y)\,S(\bold X,\bold w).
\hskip -2em
\tag9.3
$$
In the space $V$ of dimension $n\geqslant 2$ for any vector $\bold X\neq 
0$ one can find the vector $\bold Y$ noncollinear to $\bold X$. For these
two vectors there exists a linear functional such that
$$
\xalignat 2
&\alpha(\bold X)=0,&&\alpha(\bold Y)=1.\hskip -2em
\tag9.4
\endxalignat
$$
Substituting \thetag{9.4} into \thetag{9.3} we obtain $S(\bold X,\bold w)=0$
for any vector $\bold X\in V$ and for any vector $\bold w\in W=\Ker\hat R$. 
This completes the proof of lemma~9.1.
\qed\enddemo
    The condition $n\geqslant 2$ used in proof of lemma~9.1 follows
from the inequalities $1\leqslant m\leqslant n-1$ determining $m$-th case 
of intermediate degeneration. If $n=2$, then skew-symmetric bilinear
form with nonzero kernel is identically zero: $\bold S=2\,\tilde\bold R=0$.
\proclaim{Lemma 9.2} Let \thetag{1.1} be a system of equations belonging 
to $m$-th case of intermediate degeneration such that its algebra of point 
symmetries has maximal dimension \thetag{9.1}. Under these assumptions
if $m=1$ or $m\geqslant 3$, then $\bold S=2\,\tilde\bold R=0$.
\endproclaim
\demo{Proof} Let's use the result of previous lemma~9.1. Due to the 
inclusion $W\subset\Ker S$ tensor field $\bold S$ induces a skew-symmetric 
bilinear form $\tilde S$ on the factorspace $V/W$. It is defined by means
of the relationship
$$
\tilde S(\hat\bold X,\hat\bold Y)=S(\bold X,\bold Y),\hskip -2em
\tag9.5
$$
where $\hat\bold X=\Cl_{\sssize W}(\bold X)$ and $\hat\bold Y=\Cl_{
\sssize W}(\bold Y)$. For $m=1$ we have $\dim(V/W)=1$ and remember
that skew-symmetric bilinear form in one-dimensional space is 
identically zero. Therefore $\tilde S=0$. Due to \thetag{9.5} this
implies $\bold S=2\,\tilde\bold R=0$.\par
     Now let $m$ be greater than 1. The equality \thetag{9.1} for the 
dimension of the algebra of point symmetries $\Cal G$ means that the 
relationship
$$
\pagebreak
\tilde S(\hat\bold F\hat\bold X,\hat\bold Y)=-\tilde S(\hat\bold X,
\hat\bold F\hat\bold Y)\hskip -2em
\tag9.6
$$
should be fulfilled for arbitrary operator $\hat\bold F\in\End(V/W)$ 
satisfying the equality $\hat R(\hat\bold F\hat\bold X,\hat\bold Y)=
-\hat R(\hat\bold X,\hat\bold F\hat\bold Y)$. By applying the 
complexification procedure to the factorspace $V/W$, if necessary, 
the matrix of nondegenerate symmetric bilinear form $\hat R$ from 
\thetag{8.7} can be brought to the unit matrix at the expense of 
proper choice of base $\bold e_1,\,\ldots,\,\bold e_m$ in factorspace
$V/W$: 
$$
\hat R(\bold e_i,\bold e_j)=\cases 1 &\text{for \ }i=j,\\
0 &\text{for \ }i\neq j.\endcases\hskip -2em
\tag9.7
$$
When $m\geqslant 3$ we define an operator $\hat\bold F$ by prescribing
its action upon the base vectors $\bold e_1,\,\ldots,\,\bold e_m$ of
the above base in factorspace $V/W$:
$$
\hat\bold F(\bold e_i)=\cases\hphantom{-}\bold e_q &\text{for \ }i=k,\\
-\bold e_k &\text{for \ }i=q,\\ \hphantom{-}0 &\text{for \ }i\neq k
\text{\ \ and \ }i\neq q.\endcases\hskip -2em
\tag9.8
$$
From \thetag{9.7} and \thetag{9.8} for operator $\hat\bold F$ we derive
$\hat R(\hat\bold F\hat\bold X,\hat\bold Y)=-\hat R(\hat\bold X,\hat
\bold F\hat\bold Y)$. Hence \thetag{9.6} should be fulfilled for $\hat
\bold F$. Taking into account \thetag{9.8} from \thetag{9.6} we obtain 
$$
\tilde S_{ik}=\tilde S(\bold e_i,\bold e_k)=-\tilde S(\bold e_i,
\hat\bold F\bold e_q)=\tilde S(\hat\bold F\bold e_i,\bold e_q)=
S(0,\bold e_q)=0.
$$
Since $i$, $k$, $q$ are three arbitrary indices, the above equality
yields $\tilde S_{ik}=0$ for all nondiagonal elements in the matrix
of bilinear form $\tilde S$. Diagonal elements are zero due to skew
symmetry. Hence $\tilde S=0$ and $\bold S=2\,\tilde\bold R=0$. 
Lemma is proved.\qed\enddemo
\head
10. Structure of curvature tensor for $n\geqslant 4$.
\endhead
    According to the theorem~4.1 Lie derivative of curvature tensor
$\bold R$ along any field of point symmetry of the system of equations
\thetag{1.1} is zero. In local coordinates the equation $L_{\boldsymbol
\eta}(\bold R)=0$ is written in the following form:
$$
\sum^n_{k=1}\left(\eta^k\,\frac{\partial R^i_{jrs}}{\partial y^k}
-R^k_{jrs}\,\frac{\partial\eta^i}{\partial y^k}+R^i_{krs}\,
\frac{\partial\eta^k}{\partial y^j}+R^i_{jks}\,\frac{\partial\eta^k}
{\partial y^r}+R^i_{jrk}\,\frac{\partial\eta^k}{\partial y^s}\right)=0.
\hskip -2em
\tag10.1
$$
Let's take and fix some point $p_0$ on $M$. For the field of point 
symmetry $\boldsymbol\eta$ from $\Cal G(p_0)$ we can transform 
\thetag{10.1} to the form analogous to \thetag{8.3}:
$$
\bold R(\bold F\bold X,\bold Y)\bold Z+\bold R(\bold X,\bold F\bold Y)
\bold Z+\bold R(\bold X,\bold Y)\bold F\bold Z=\bold F\bold R(\bold X,
\bold Y)\bold Z.\hskip -3em
\tag10.2
$$
Here $\bold X$, $\bold Y$, $\bold Z$ are three arbitrary vectors from
tangent space $V=T_{p_0}(M)$. Let \thetag{1.1} be a system of equations
belonging to $m$-th case of intermediate degeneration. Then \thetag{9.1}
means that the relationship \thetag{10.2} is fulfilled for all those
operators $\bold F$, for which \thetag{8.3} holds.\par
     In $m$-th case of intermediate degeneration the kernel of bilinear
form $\hat R$ in \thetag{8.3} is nonzero since $1\leqslant m\leqslant n-1$.
\pagebreak
Choosing some nonzero vector $\bold w\in W=\Ker\hat R$ we construct an
operator $\bold F=\bold w\otimes\alpha$, where $\alpha$ is an arbitrary
linear functional in $V$. Such operator satisfies the equality \thetag{8.3}. 
By substituting it into \thetag{10.2} we get
$$
\aligned
\alpha(\bold X)&\cdot\bold R(\bold w,\bold Y)\bold Z+\alpha(\bold Y)
\cdot\bold R(\bold X,\bold w)\bold Z\,+\\
\vspace{1ex}
&+\,\alpha(\bold Z)\cdot\bold R(\bold X,\bold Y)\bold w=\alpha(\bold R
(\bold X,\bold Y)\bold Z)\cdot\bold w.
\endaligned\hskip -3em
\tag10.3
$$
Formula \thetag{10.3} is a key to further analysis of the structure
of curvature tensor $\bold R$.
\proclaim{Theorem 10.1} Let \thetag{1.1} be a system of equations 
belonging to $m$-th case of intermediate degeneration such that its 
algebra of point symmetries has maximal dimension \thetag{9.1}. 
Under these assumptions for $n\geqslant 4$ there is a formula
$$
\bold R(\bold X,\bold Y)\bold Z=\sigma(\bold X,\bold Y)\cdot\bold Z+
\beta(\bold Y,\bold Z)\cdot\bold X-\beta(\bold X,\bold Z)\cdot\bold Y
\hskip -2em
\tag10.4
$$
expressing curvature tensor through two tensor fields $\boldsymbol\sigma$ 
and $\boldsymbol\beta$ of type $(0,2)$.
\endproclaim
     For three arbitrary vectors $\bold X$, $\bold Y$, $\bold Z$ in
\thetag{10.3} let's consider their linear span $U=\langle\bold X,\bold Y,
\bold Z\rangle$. Due to $\dim V=n\geqslant 4$ subspace $U$ do not 
coincide 
with $V$. Denote by $U^\perp$ the set of linear functionals $\alpha$ such 
that $\alpha(\bold X)=0$, $\alpha(\bold Y)=0$, and $\alpha(\bold Z)=0$ 
simultaneously. Then $U^\perp$ is a subspace in dual space $V^*$ such 
that $\dim U^\perp=n-\dim U$ (see for instance \cite{18}) and
$$
U=\{\bold u\in V:\ \ \alpha(\bold u)=0\text{\ \ for all \ }
\alpha\in U^\perp\}.\hskip -3em
\tag10.5
$$
Recall that $\alpha$ in \thetag{10.3} is an arbitrary linear functional.
Substituting various linear functionals $\alpha\in U^\perp$ into 
\thetag{10.3} we find that $\alpha(\bold R(\bold X,\bold Y)\bold Z)=0$
for them. Due to \thetag{10.5} this means that for any three vectors
$\bold X$, $\bold Y$, $\bold Z$ vector $\bold R(\bold X,\bold Y)\bold Z$
is in their linear span. Let's express this circumstance as
$$
\bold R(\bold X,\bold Y)\bold Z=\beta\cdot\bold X+\gamma\cdot\bold Y+
\sigma\cdot\bold Z.\hskip -3em
\tag10.6
$$
Now denote by $U=\langle\bold X,\bold Y\rangle$ linear span of
two vectors $\bold X$ and $\bold Y$. If $\bold Z\notin U$, then 
coefficient $\sigma$ is uniquely defined by the expansion \thetag{10.6}.
Thus we have a function $\sigma=\sigma(\bold X,\bold Y,\bold Z)$ defined
for triples of vectors such that $\bold Z\notin\langle\bold X,\bold Y
\rangle$. Further proof of theorem~10.1 breaks into series of lemmas. 
\proclaim{Lemma 10.1} For $\dim V\geqslant 4$ the function $\sigma=
\sigma(\bold X,\bold Y,\bold Z)$ doesn't depend on $\bold Z$.
\endproclaim
\demo{Proof} Let's retain the notation $U=\langle\bold X,\bold Y\rangle$
for linear span of the vectors $\bold X$ and $\bold Y$ and consider a
factorspace $V/U$. Let $\bold Z_1$ and $\bold Z_2$ be two arbitrary 
vectors such that their cosets relative to subspace $U$ are linearly 
independent. The existence of such vectors $\bold Z_1$ and $\bold Z_2$
follows from the estimate 
$$
\dim(V/U)\geqslant 4-2=2.
$$
For these two vectors we have \pagebreak $\bold Z_1\notin U$ and $\bold Z_2
\notin U$ so that if $\bold Z_3=\bold Z_1+\bold Z_2$, then $\bold Z_3\notin 
U$. Let's write the equation  \thetag{10.6}  for  each  of  three 
triples of 
vectors:
$$
\aligned
\bold R(\bold X,\bold Y)\bold Z_1=\beta_1\cdot\bold X+\gamma_1\cdot\bold Y+
\sigma_1\cdot\bold Z_1,\\
\vspace{1ex}
\bold R(\bold X,\bold Y)\bold Z_2=\beta_2\cdot\bold X+\gamma_2\cdot\bold Y+
\sigma_2\cdot\bold Z_2,\\
\vspace{1ex}
\bold R(\bold X,\bold Y)\bold Z_3=\beta_3\cdot\bold X+\gamma_3\cdot\bold Y+
\sigma_3\cdot\bold Z_3.
\endaligned\hskip -3em
\tag10.7
$$
Let's add first two equalities \thetag{10.7} and subtract the third one.
Then factorize the obtained equality with respect to the subspace $U$.
As a result we have 
$$
(\sigma_1-\sigma_3)\cdot\Cl_{\sssize U}(\bold Z_1)+
(\sigma_2-\sigma_3)\cdot\Cl_{\sssize U}(\bold Z_2)=0.
$$
Since cosets $\Cl_{\sssize U}(\bold Z_1)$ and $\Cl_{\sssize U}(\bold 
Z_2)$ are linearly independent, from the above equality we obtain 
$\sigma_1=\sigma_3$ and $\sigma_2=\sigma_3$, where $\sigma_1=\sigma(
\bold X,\bold Y,\bold Z_1)$ and $\sigma_2=\sigma(\bold X,\bold Y,\bold 
Z_2)$. Thus we have proved the required result $\sigma(\bold X,\bold Y,
\bold Z_1)=\sigma(\bold X,\bold Y,\bold Z_2)$ for the vectors $\bold Z_1$ 
and $\bold Z_2$ whose cosets are linearly independent.\par
     Now suppose that cosets $\Cl_{\sssize U}(\bold Z_1)$ and $\Cl_{
\sssize U}(\bold Z_1)$ are linearly dependent, but are nonzero. Then
due to $\dim(V/U)\geqslant 2$ we can find a vector $\bold Z_4$ coset of
which is not collinear to cosets $\Cl_{\sssize U}(\bold Z_1)$ and $\Cl_{
\sssize U}(\bold Z_2)$, and we can apply previous result:
$$
\sigma(\bold X,\bold Y,\bold Z_1)=\sigma(\bold X,\bold Y,\bold Z_4)=
\sigma(\bold X,\bold Y,\bold Z_2).
$$
Cases when $\Cl_{\sssize U}(\bold Z_1)=0$ or $\Cl_{\sssize U}(\bold 
Z_2)=0$ are not considered since in these cases the equality 
\thetag{10.6} do not define both quantities $\sigma(\bold X,\bold Y,
\bold Z_1)$ and $\sigma(\bold X,\bold Y,\bold Z_2)$.
\qed\enddemo
     Lemma~10.1 shows that for $n\geqslant 4$ the equality \thetag{10.6}
defines the function $\sigma=\sigma(\bold X,\bold Y)$ which do not depend
on third vector $\bold Z$. This circumstance gives us the opportunity to 
expand \thetag{10.6} for the case when $\bold Z$ belongs to $U=\langle
\bold X,\bold Y\rangle$.
\proclaim{Lemma 10.2} For $\dim V\geqslant 4$ the function $\sigma(\bold X,
\bold Y)$ defined by the relationship \thetag{10.6} is linear in its first 
argument $\bold X$.
\endproclaim
\demo{Proof} Let $\bold X_1+\bold X_2=\bold X_3$. Denote by $\bold U=\langle
\bold X_1,\bold X_2,\bold X_3,\bold Y\rangle$ linear span of four vectors
$\bold X_1$, $\bold X_2$, $\bold X_3$, and $\bold Y$. Its obvious that
$\dim U\leqslant 3$. Due to $\dim V\geqslant 4$ we can find a vector 
$\bold Z\notin U$. Hence $\Cl_{\sssize U}(\bold Z)\neq 0$, and we can
write \thetag{10.6} three times:
$$
\aligned
\bold R(\bold X_1,\bold Y)\bold Z=\beta_1\cdot\bold X_1+\gamma_1\cdot
\bold Y+\sigma(\bold X_1,\bold Y)\cdot\bold Z,\\
\vspace{1ex}
\bold R(\bold X_2,\bold Y)\bold Z=\beta_2\cdot\bold X_2+\gamma_2\cdot
\bold Y+\sigma(\bold X_2,\bold Y)\cdot\bold Z,\\
\vspace{1ex}
\bold R(\bold X_3,\bold Y)\bold Z=\beta_3\cdot\bold X_3+\gamma_3\cdot
\bold Y+\sigma(\bold X_2,\bold Y)\cdot\bold Z.
\endaligned\hskip -3em
\tag10.8
$$
Let's add first two equalities \thetag{10.8} and subtract the third one.
Then factorize with respect to the subspace $U$. As a result we have
$$
\pagebreak
\sigma(\bold X_1+\bold X_2,\bold Y)=\sigma(\bold X_1,\bold Y)
+\sigma(\bold X_2,\bold Y).\hskip -3em
\tag10.9
$$\par
    Now suppose that we are given a real number $\alpha\in\Bbb R$ and
two vectors $\bold X_1$ and $\bold Y$. Suppose $\bold X_2=\alpha\cdot
\bold X_1$ and consider linear span $U=\langle\bold X_1,\bold X_2,\bold 
Y\rangle$ for the vectors $\bold X_1$, $\bold X_2$, and $\bold Y$. Due
to inequalities $\dim V\geqslant 4$ and $\dim U\leqslant 2$ we can find 
a vector $\bold Z\notin U$. Then write the equality \thetag{10.6} twice:
$$
\aligned
\bold R(\bold X_1,\bold Y)\bold Z=\beta_1\cdot\bold X_1+\gamma_2\cdot
\bold Y+\sigma(\bold X_1,\bold Y)\cdot\bold Z,\\
\vspace{1ex}
\bold R(\bold X_2,\bold Y)\bold Z=\beta_2\cdot\bold X_2+\gamma_2\cdot
\bold Y+\sigma(\bold X_2,\bold Y)\cdot\bold Z.
\endaligned\hskip -3em
\tag10.10
$$
Multiply first equality \thetag{10.10} by $\alpha$ and subtract second 
one from it. After factorizing with respect to the subspace $U$ we 
finally get
$$
\sigma(\alpha\cdot\bold X_1,\bold Y)=\alpha\,\sigma(\bold X_1,\bold Y). 
$$
This equality together with \thetag{10.9} proves linearity of $\sigma(
\bold X,\bold Y)$ in its first argument. Lemma is proved.
\qed\enddemo
\proclaim{Lemma 10.3} For $\dim V\geqslant 4$ the function $\sigma(\bold X,
\bold Y)$ defined by the relationship \thetag{10.6} is linear in its second
argument $\bold Y$.
\endproclaim
    Proof of the lemma~10.3 is quite similar to the proof of previous 
lemma~10.2 so that we can omit it. In whole, lemmas~10.1, 10.2, and 10.3
mean that the relationship \thetag{10.6} defines a bilinear form
$\sigma(\bold X,\bold Y)$ and corresponding tensor $\boldsymbol\sigma$
of type $(0,2)$. Therefore we rewrite \thetag{10.6} as
$$
\bold R(\bold X,\bold Y)\bold Z-\sigma(\bold X,\bold Y)\cdot\bold Z
=\beta\cdot\bold X+\gamma\cdot\bold Y.\hskip -3em
\tag10.11
$$
Note that for to state and prove lemmas~10.1, 10.2, and 10.3 we used
only the linearity of left hand side of \thetag{10.6} with respect
to $\bold X$, $\bold Y$, and $\bold Z$. Left had side of \thetag{10.11}
is also linear in $\bold X$, $\bold Y$, and $\bold Z$. Repeating
previous arguments we obtain two bilinear forms $\beta(\bold Y,\bold Z)$ 
and $\gamma(\bold X,\bold Z)$ such that \thetag{10.11} can be written as
$$
\bold R(\bold X,\bold Y)\bold Z=\sigma(\bold X,\bold Y)\cdot\bold Z+
\beta(\bold Y,\bold Z)\cdot\bold X+\gamma(\bold X,\bold Z)\cdot\bold Y.
\hskip -3em
\tag10.12
$$
Left hand side of \thetag{10.12} is skew-symmetric respective to $\bold X$
and $\bold Y$. This yields $\boldsymbol\gamma=-\boldsymbol\beta$. Therefore
we can write \thetag{10.12} in form of \thetag{10.4}, thus completing 
proof of the theorem~10.1. Furthermore skew symmetry of $\bold R(\bold 
X,\bold Y)$ in $\bold X$ and $\bold Y$ implies skew symmetry of bilinear 
form $\sigma(\bold X,\bold Y)$.\par
     Formula \thetag{10.4} gives us the opportunity to express the
components of curvature tensor through the components of two tensor 
fields $\boldsymbol\beta$ and $\boldsymbol\sigma$:
$$
R^k_{sij}=\sigma_{ij}\,\delta^k_s+\beta_{js}\,\delta^k_i-\beta_{is}\,
\delta^k_j.\hskip -3em
\tag10.13
$$
Now if we make contractions in \thetag{10.13} according to the formulas
\thetag{4.11}, we obtain the following expressions for the components of
tensor $\bold S$ and Ricci tensor $\bold R$:
$$
\pagebreak
\xalignat 2
&S_{ij}=n\,\sigma_{ij}-(\beta_{ij}-\beta_{ji}),
&&R_{ij}=\sigma_{ij}+(n-1)\,\beta_{ji}.
\endxalignat
$$
Let $\hat{\boldsymbol\beta}$ and $\tilde{\boldsymbol\beta}$ be symmetric
and skew-symmetric parts of tensor $\boldsymbol\beta$. Then
$$
\xalignat 3
&S_{ij}=n\,\sigma_{ij}-2\,\tilde\beta_{ij},
&&\hat R_{ij}=(n-1)\,\hat\beta_{ij},
&&\tilde R_{ij}=\sigma_{ij}-(n-1)\,\tilde\beta_{ij}.
\endxalignat
$$
If we take into account the relationship $\bold S=2\,\tilde\bold R$,
we can express tensor fields $\boldsymbol\sigma$, $\hat{\boldsymbol\beta}$,
and $\tilde{\boldsymbol\beta}$ through symmetric and skew-symmetric
parts of Ricci tensor:
$$
\xalignat 3
&\hat{\boldsymbol\beta}=\frac{\hat\bold R}{n-1},
&&\tilde{\boldsymbol\beta}=-\frac{\tilde\bold R}{n+1},
&&\boldsymbol\sigma=\frac{2\,\tilde\bold R}{n+1}.
\hskip -3em
\tag10.14
\endxalignat
$$
From \thetag{10.13} and \thetag{10.14} we derive formula for components 
of curvature tensor:
$$
R^k_{sij}=\frac{2\,\tilde R_{ij}\,\delta^k_s-\tilde R_{js}\,\delta^k_i
+\tilde R_{is}\,\delta^k_j}{n+1}+\frac{\hat R_{js}\,\delta^k_i
-\hat R_{is}\,\delta^k_j}{n-1}.\hskip -3em
\tag10.15
$$
\proclaim{Theorem 10.2} Let \thetag{1.1} be a system of equations belonging 
to $m$-th case of intermediate degeneration such that its algebra of point 
symmetries has maximal dimension \thetag{9.1}. Under these assumptions
if $n\geqslant 4$, then appropriate curvature tensor is expressed through 
Ricci tensor according to the formula \thetag{10.15}.
\endproclaim
\head
11. Structure of curvature tensor for $n=3$.
\endhead
    Now we continue considering cases of intermediate degeneration
when $1\leqslant m\leqslant n-1$. Dimension $n=3$ restricts our 
possibilities to maneuver in proving propositions like lemmas~10.1 
and 10.2. In order to avoid these restrictions let's substitute
$\bold Z=\bold X$ into the formula \thetag{10.3}. This yields
$$
\aligned
\alpha(\bold X)&\cdot\bold R(\bold w,\bold Y)\bold X+\alpha(\bold Y)
\cdot\bold R(\bold X,\bold w)\bold X\,+\\
\vspace{1ex}
&+\,\alpha(\bold X)\cdot\bold R(\bold X,\bold Y)\bold w=\alpha(\bold R
(\bold X,\bold Y)\bold X)\cdot\bold w.
\endaligned\hskip -3em
\tag11.1
$$
\proclaim{Theorem 11.1} Let \thetag{1.1} be a system of equations belonging 
to $m$-th case of intermediate degeneration such that its algebra of point 
symmetries has maximal dimension \thetag{9.1}. Under these assumptions
if $n=3$, then there is a relationship
$$
\bold R(\bold X,\bold Y)\bold X=\theta(\bold X,\bold Y)\cdot\bold X
-\beta(\bold X,\bold X)\cdot\bold Y\hskip -3em
\tag11.2
$$
binding curvature tensor with two tensor fields $\boldsymbol\beta$ and 
$\boldsymbol\theta$ of type $(0,2)$.
\endproclaim
     For two arbitrary vectors $\bold X$ and $\bold Y$ we consider
their linear span $U=\langle\bold X,\bold Y\rangle$. Due to $\dim V
=n=3$ subspace $U$ do not coincide with $V$. Denote by $U^\perp$ the 
set of linear functionals $\alpha$ such that $\alpha(\bold X)=0$ and
$\alpha(\bold Y)=0$ simultaneously. Then $U^\perp$ is a subspace in 
dual space $V^*$ such that 
$$
U=\{\bold u\in V:\ \ \alpha(\bold u)=0\text{\ \ for all \ }
\alpha\in U^\perp\}.\hskip -3em
\tag11.3
$$
Recall that $\alpha$ in \thetag{11.1} is an arbitrary linear functional.
\pagebreak Substituting various linear functionals $\alpha\in U^\perp$ 
into \thetag{11.1} we find that $\alpha(\bold R(\bold X,\bold Y)\bold X)=0$
for them. Due to \thetag{11.3} this means that for any two vectors
$\bold X$ and $\bold Y$ vector $\bold R(\bold X,\bold Y)\bold X$
is in their linear span. Let's express this circumstance as
$$
\bold R(\bold X,\bold Y)\bold X=\theta\cdot\bold X
-\beta\cdot\bold Y.\hskip -3em
\tag11.4
$$
Now denote by $U=\langle\bold X\rangle$ linear span of the vector 
$\bold X$. If $\bold Y\notin U$, then coefficient $\beta$ is uniquely 
defined by the expansion \thetag{11.4}. Thus we have a function 
$\beta=\beta(\bold X,\bold Y)$ defined for pairs of vectors such that 
$\bold Y\notin\langle\bold X\rangle$. Further proof of theorem~11.1 
breaks into series of lemmas. 
\proclaim{Lemma 11.1} For $n=3$ the function $\beta=\beta(\bold X,
\bold Y)$ doesn't depend on $\bold Y$.
\endproclaim
\demo{Proof} Let's retain the notation $U=\langle\bold X,\bold Y\rangle$
for linear span of the vector $\bold X$ and consider a factorspace $V/U$. 
Let $\bold Y_1$ and $\bold Y_2$ be two arbitrary vectors such that their 
cosets relative to subspace $U$  are  linearly  independent.  The 
existence 
of such vectors $\bold Y_1$ and $\bold Y_2$ follows from the estimate 
$$
\dim(V/U)\geqslant 3-1=2.
$$
For these two vectors we have $\bold Y_1\notin U$ and $\bold Y_2\notin U$ 
so that if $\bold Y_3=\bold Y_1+\bold Y_2$, then $\bold Y_3\notin 
U$. Let's write the equation \thetag{11.4} for each of three pairs of 
vectors:
$$
\aligned
\bold R(\bold X,\bold Y_1)\bold X=\theta_1\cdot\bold X-\beta_1\cdot
\bold Y_1,\\
\vspace{1ex}
\bold R(\bold X,\bold Y_2)\bold X=\theta_2\cdot\bold X-\beta_2\cdot
\bold Y_2,\\
\vspace{1ex}
\bold R(\bold X,\bold Y_3)\bold X=\theta_3\cdot\bold X-\beta_3\cdot
\bold Y_3.
\endaligned\hskip -3em
\tag11.5
$$
Let's add first two equalities \thetag{11.5} and subtract the third one.
Then factorize the obtained equality with respect to the subspace $U$.
As a result we have 
$$
(\beta_1-\beta_3)\cdot\Cl_{\sssize U}(\bold Y_1)+
(\beta_2-\beta_3)\cdot\Cl_{\sssize U}(\bold Y_2)=0.
$$
Since cosets $\Cl_{\sssize U}(\bold Y_1)$ and $\Cl_{\sssize U}(\bold 
Y_2)$ are linearly independent, from the above equality we obtain 
$\beta_1=\beta_3$ and $\beta_2=\beta_3$, where $\beta_1=\beta(\bold X,
\bold Y_1)$ and $\beta_2=\beta(\bold X,\bold Y_2)$. 
Thus we have proved the required result $\beta(\bold X,\bold Y_1)=
\beta(\bold X,\bold Y_2)$ for the vectors $\bold Y_1$ and $\bold Y_2$,
whose cosets are linearly independent.\par
     Now suppose that cosets $\Cl_{\sssize U}(\bold Y_1)$ and $\Cl_{
\sssize U}(\bold Y_1)$ are linearly dependent, but are nonzero. Then
due to $\dim(V/U)\geqslant 2$ we can find a vector $\bold Y_4$ coset of
which is not collinear to cosets $\Cl_{\sssize U}(\bold Y_1)$ and $\Cl_{
\sssize U}(\bold Y_2)$, and we can apply previous result:
$$
\beta(\bold X,\bold Y_1)=\beta(\bold X,\bold Y_4)=\beta(\bold X,\bold Y_2).
$$
Cases when $\Cl_{\sssize U}(\bold Y_1)=0$ or $\Cl_{\sssize U}(\bold Y_2)=0$ 
are not considered since in these cases the equality \thetag{11.4} do not 
define both quantities $\beta(\bold X,\bold Y_1)$ and $\beta(\bold X,
\bold Y_2)$.\qed\enddemo
     Lemma~11.1 shows that for the dimension $n=3$ the equality 
\thetag{11.4} defines the function $\beta=\beta(\bold X,\bold Y)$, 
\pagebreak which do not depend on second vector $\bold Y$. This 
circumstance gives us the opportunity to expand \thetag{11.4} for 
the case when $\bold Y$ belongs to linear span of vector $\bold X$.
\par
     Remember that numerical function $\beta(\bold X)$ of one vectorial
argument is called quadratic form if it is obtained from some bilinear
form $\gamma(\bold X,\bold Z)$ by substitution $\bold Z=\bold X$, i\.~e\.
$\beta(\bold X)=\gamma(\bold X,\bold X)$. Without loss of generality
bilinear form $\gamma$ can be assumed to be symmetric (see more details
in \cite{18}).
\proclaim{Lemma 11.2} Numerical function of one vectorial argument
$\beta(\bold X)$ is a quadratic form if and only if it satisfies the
relationship
$$
\beta(\bold X+\alpha\cdot\bold Z)+\alpha\,\beta(\bold X-\bold Z)=
(1+\alpha)\,\beta(\bold X)+(\alpha+\alpha^2)\,\beta(\bold Z),
\hskip -3em
\tag11.6
$$
where $\bold X$ and $\bold Z$ are two arbitrary vectors and $\alpha$
is an arbitrary number.
\endproclaim
\demo{Proof} Direct statement of lemma is obvious. Indeed, if for some 
bilinear form $\gamma$ we substitute $\beta(\bold X)=\gamma(\bold X,
\bold X)$ into the relationship \thetag{11.6}, we find that it is 
fulfilled identically.\par
    Now let's prove converse assertion. Suppose \thetag{11.6} to
be fulfilled for the function $\beta(\bold X)$. Substituting $\bold X
=\bold Z=0$ and $\alpha=1$ into \thetag{11.6} we get 
$$
\beta(0)=0.\hskip -3em
\tag11.7
$$
As a second step we substitute $\bold Z=\bold X$ into \thetag{11.6}.
If we take into account \thetag{11.7}, we obtain $\beta((1+\alpha)\cdot
\bold X)=(1+\alpha)^2\,\beta(\bold X)$. This equality can be simplified 
by substituting $\alpha-1$ for $\alpha$. It takes the following form:
$$
\beta(\alpha\cdot\bold X)=\alpha^2\,\beta(\bold X).\hskip -3em
\tag11.8
$$
So $\beta$ is a homogeneous function of degree $2$. Now we use 
$\beta(\bold X)$ in order to define the function $\gamma(\bold X,
\bold Z)$ of two vectorial arguments:
$$
\gamma(\bold X,\bold Z)=\frac{\beta(\bold X+\bold Z)-\beta(\bold X
-\bold Z)}{4}.\hskip -3em
\tag11.9
$$\par
     Due to \thetag{11.8} the function $\gamma(\bold X,\bold Z)$ in 
\thetag{11.9} is symmetric: $\gamma(\bold X,\bold Z)=\gamma(\bold Z,
\bold X)$. Let's prove that it is linear respective to its second 
argument. We write the relationship \thetag{11.6} substituting
$\alpha$ for $-\alpha$ in it:
$$
\beta(\bold X-\alpha\cdot\bold Z)-\alpha\,\beta(\bold X-\bold Z)=
(1-\alpha)\,\beta(\bold X)+(\alpha^2-\alpha)\,\beta(\bold Z).
\hskip -3em
\tag11.10
$$
Let's subtract \thetag{11.10} from the initial relationship \thetag{11.6}.
This yields
$$
4\,\gamma(\bold X,\alpha\cdot\bold Z)=\beta(\bold X+\alpha\cdot\bold Z)
-\beta(\bold X-\alpha\cdot\bold Z)=2\,\alpha\,(\beta(\bold X)+\beta(
\bold Z)-\beta(\bold X-\bold Z)).
$$
For $\alpha=1$ this relationship takes the form
$$
\pagebreak
4\,\gamma(\bold X,\bold Z)=\beta(\bold X+\bold Z)-\beta(\bold X-\bold Z)
=2\,(\beta(\bold X)+\beta(\bold Z)-\beta(\bold X-\bold Z)).
$$
Comparing two above relationships we find that 
$$
\gamma(\bold X,\alpha\cdot\bold Z)=\alpha\,\gamma(\bold X,\bold Z).
\hskip -3em
\tag11.11
$$
For $\alpha=1$ the relationship \thetag{11.6} can be written as follows:
$$
2\,\beta(\bold X)+2\,\beta(\bold Z)=\beta(\bold X+\bold Z)
+\beta(\bold X-\bold Z).\hskip -3em
\tag11.12
$$
Let's substitute $\bold X$ by $\bold X+\bold Z_1$ and $\bold Z$ by
$\bold X+\bold Z_2$ in \thetag{11.12}. As a result we obtain
$$
2\,\beta(\bold X+\bold Z_1)+2\,\beta(\bold X+\bold Z_2)=\beta(2\cdot
\bold X+\bold Z_1+\bold Z_2)+\beta(\bold Z_1-\bold Z_2).\hskip -3em
\tag11.13
$$
Further in formula \thetag{11.12} we substitute $\bold X$ by $\bold X
-\bold Z_1$ and $\bold Z$ by $\bold X-\bold Z_2$:
$$
2\,\beta(\bold X-\bold Z_1)+2\,\beta(\bold X-\bold Z_2)=\beta(2\cdot\bold X
-\bold Z_1-\bold Z_2)+\beta(\bold Z_1-\bold Z_2).\hskip -3em
\tag11.14
$$
Then we subtract \thetag{11.14} from \thetag{11.13} and compare the
expressions in both sides of the obtained formula with \thetag{11.9}.
As a result of this comparison we find
$$
8\,\gamma(\bold X,\bold Z_1)+8\,\gamma(\bold X,\bold Z_2)=
4\,\gamma(2\cdot\bold X,\bold Z_1+\bold Z_2).
$$
If we take into account \thetag{11.11} and the symmetry of the function
$\gamma(\bold X,\bold Z)$, we can bring this relationship to the following
form:
$$
\gamma(\bold X,\bold Z_1)+\gamma(\bold X,\bold Z_2)=\gamma(\bold X,
\bold Z_1+\bold Z_2).\hskip -3em
\tag11.15
$$
The relationships \thetag{11.11} and \thetag{11.15} in the aggregate
mean the linearity of $\gamma(\bold X,\bold Z)$ respective to $\bold Z$.
Linearity of $\gamma(\bold X,\bold Z)$ respective to $\bold X$ then
follows by the symmetry $\gamma(\bold X,\bold Z)=\gamma(\bold Z,
\bold X)$. Thus we constructed symmetric bilinear form $\gamma(\bold X,
\bold Z)$. Due to \thetag{11.7} and \thetag{11.8} if we substitute
$\bold Z=\bold X$ in $\gamma(\bold X,\bold Z)$, we get $\beta(\bold X)
=\gamma(\bold X,\bold X)$. Proof is complete.\qed\enddemo
     Now let's return to the formula \thetag{11.4}. Due to the lemma~11.1
proved above this formula can be written as:
$$
\bold R(\bold X,\bold Y)\bold X=\theta\cdot\bold X
-\beta(\bold X)\cdot\bold Y.\hskip -3em
\tag11.16
$$
Left hand side of \thetag{11.16} satisfies the relationship
$$
\gathered
\bold R(\bold X+\alpha\cdot\bold Z,\bold Y)(\bold X+\alpha\cdot
\bold Z)+\alpha\,\bold R(\bold X+\bold Z,\bold Y)(\bold X+\bold Z)=\\
\vspace{1ex}
=(1+\alpha)\,\bold R(\bold X,\bold Y)\bold X+(\alpha+\alpha^2)\,
\bold R(\bold Z,\bold Y)\bold Z,
\endgathered\hskip -3em
\tag11.17
$$
which can be checked by direct calculations. For the case $\dim V=3$
we can choose vector $\bold Y$ such that it doesn't belong to the
linear span of $\bold X$ and $\bold Z$. Therefore, if we transform 
each summand in \thetag{11.17} by means of \thetag{11.16}, we obtain 
vectorial equality which leads to the relationship \thetag{11.6} for 
the function $\beta(\bold X)$.
\proclaim{Lemma 11.3} For the dimension $n=3$ function $\beta(\bold X)$
in \thetag{11.16} is a quadratic form, which is defined by some bilinear 
form.
\endproclaim
    Lemma~11.3 is an immediate consequence of the lemma~11.2 due to above
arguments. We express the result of this lemma as follows:
$$
\bold R(\bold X,\bold Y)\bold X+\beta(\bold X,\bold X)\cdot\bold Y
=\theta\cdot\bold X.\hskip -3em
\tag11.18
$$
\proclaim{Lemma 11.4} For the dimension $n=3$ the coefficient 
$\theta$ in \thetag{11.18} is a bilinear function of two vectorial 
arguments $\bold X$ and $\bold Y$.
\endproclaim
\demo{Proof} In order to prove the linearity o $\theta(\bold X,\bold Y)$
with respect to $\bold Y$ we should note that for any fixed vector $\bold X
\neq 0$ in the space $V$ of dimension $n=3$ one can find linear functional 
$\gamma$ such that $\gamma(\bold X)=1$. Applying this functional to both 
sides of the equality \thetag{11.18} we obtain
$$
\theta(\bold X,\bold Y)=\gamma(\bold R(\bold X,\bold Y)\bold X)+
\beta(\bold X,\bold X)\,\gamma(\bold Y).
$$
Linearity of $\theta(\bold X,\bold Y)$ respective to $\bold Y$ now
is a consequence of linearity of right hand side of the above equality 
respective to $\bold Y$ for fixed $\gamma$ and $\bold X$.\par
    Let's prove linearity of $\theta(\bold X,\bold Y)$ respective to 
$\bold X$. Suppose $\bold X\neq 0$. If we substitute $\alpha\cdot
\bold X$ for $\bold X$ into the formula \thetag{11.18}, we get 
$$
\alpha^2\cdot\bold R(\bold X,\bold Y)\bold X+\alpha^2\,\beta(\bold X,
\bold X)\cdot\bold Y=\alpha\,\theta(\alpha\cdot\bold X,\bold Y)\cdot
\bold X. 
$$
Comparing this with the equality $\bold R(\bold X,\bold Y)\bold X
+\beta(\bold X,\bold X)\cdot\bold Y=\theta(\bold X,\bold Y)\cdot
\bold X$ we find
$$
\theta(\alpha\cdot\bold X,\bold Y)=\alpha\,\theta(\bold X,\bold Y).
\hskip -3em
\tag11.19
$$
For $\bold X=0$ the quantity $\theta$ is not correctly defined from
\thetag{11.18}. The relationship \thetag{11.19} makes possible to
extend the function $\theta(\bold X,\bold Y)$ for $\bold X=0$ by 
the value $\theta(0,\bold Y)=0$.\par
     Let $\bold X_1$ and $\bold X_2$ be two arbitrary vectors. Denote
$\bold X_3=\bold X_1+\bold X_2$. If $\bold X_1$ and $\bold X_2$ are
collinear, then there is a vector $\bold e\neq 0$ such that
$$
\xalignat 3
&\bold X_1=\alpha_1\cdot\bold e,
&&\bold X_2=\alpha_2\cdot\bold e,
&&\bold X_3=\alpha_3\cdot\bold e,
\endxalignat
$$
where $\alpha_3=\alpha_1+\alpha_2$. By means of \thetag{11.19} we
obtain $\theta(\bold X_1,\bold Y)+\theta(\bold X_2,\bold Y)=\alpha_1\,
\theta(\bold e,\bold Y)+\alpha_2\,\theta(\bold e,\bold Y)=\alpha_3\,
\theta(\bold e,\bold Y)=\theta(\bold X_3,\bold Y)$. This means
$$
\theta(\bold X_1,\bold Y)+\theta(\bold X_2,\bold Y)=\theta(\bold X_1
+\bold X_2,\bold Y).\hskip -3em
\tag11.20
$$\par
    Now suppose that vectors $\bold X_1$ and $\bold X_2$ are not
collinear. Let $\gamma$ be an arbitrary linear functional in $V$.
Let's apply it to both sides of \thetag{11.18}. As a result we
get
$$
\pagebreak
\gamma(\bold R(\bold X,\bold Y)\bold X)+\beta(\bold X,\bold X)\,\gamma(
\bold Y)=\theta(\bold X,\bold Y)\,\gamma(\bold X).\hskip -3em
\tag11.21
$$
For the fixed vector $\bold Y$ left hand side of \thetag{11.21} is a
quadratic form respective to $\bold X$ (it is obtained from bilinear
function $\gamma(\bold R(\bold X,\bold Y)\bold Z)+\beta(\bold X,\bold 
Z)\,\gamma(\bold Y)$ by substitution $\bold Z=\bold X$). Therefore
right hand side of \thetag{11.21} is also a quadratic form respective 
to $\bold X$. According to the lemma~11.2 it satisfies the relationship
\thetag{11.6} where we can take $\alpha=1$, $\bold X=\bold X_1$, and 
$\bold Z=\bold X_2$. This yields
$$
\gathered
\theta(\bold X_1+\bold X_2,\bold Y)\,\gamma(\bold X_1+\bold X_2)+
\theta(\bold X_1-\bold X_2,\bold Y)\,\gamma(\bold X_1-\bold X_2)=\\
\vspace{1ex}
2\,\theta(\bold X_1,\bold Y)\,\gamma(\bold X_1)+2\,\theta(\bold X_2,
\bold Y)\,\gamma(\bold X_2).
\endgathered\hskip -3em
\tag11.22
$$
Relying on the linearity of functional $\gamma$ we can bring \thetag{11.22} 
to the form
$$
\gathered
(\theta(\bold X_1+\bold X_2,\bold Y)+\theta(\bold X_1-\bold X_2,\bold Y)
-2\,\theta(\bold X_1,\bold Y))\,\gamma(\bold X_1)\,+\\
\vspace{1ex}
+\,(\theta(\bold X_1+\bold X_2,\bold Y)-\theta(\bold X_1-\bold X_2,
\bold Y)-2\,\theta(\bold X_2,\bold Y))\,\gamma(\bold X_2)=0.
\endgathered\hskip -3em
\tag11.23
$$
Functional $\gamma$ in \thetag{11.23} is an arbitrary linear functional.
For the pair of noncollinear vectors $\bold X_1$ and $\bold X_2$ we 
can find pair of functionals $\gamma_2$ and $\gamma_2$ such that
$$
\xalignat 4
&\gamma_1(\bold X_1)=1,&&\gamma_1(\bold X_2)=0,
&&\gamma_2(\bold X_1)=0,&&\gamma_2(\bold X_2)=1.
\endxalignat
$$
Therefore the relationship \thetag{11.23} can be broken into two
separate parts:
$$
\aligned
&\theta(\bold X_1+\bold X_2,\bold Y)+\theta(\bold X_1-\bold X_2,\bold Y)
=2\,\theta(\bold X_1,\bold Y),\\
\vspace{1ex}
&\theta(\bold X_1+\bold X_2,\bold Y)-\theta(\bold X_1-\bold X_2,\bold Y)
=2\,\theta(\bold X_2,\bold Y).
\endaligned
$$
If we add these two equalities and divide the result by $2$, we come
to the relationship \thetag{11.20}, which appears to be fulfilled
for noncollinear vectors $\bold X_1$ and $\bold X_2$ as well. In the 
aggregate with \thetag{11.19} it means that $\theta(\bold X,\bold Y)$
is linear respective to its first argument $\bold X$. Lemma is proved.
\qed\enddemo
     Lemma~11.4 completes the proof of the theorem~11.1. We shall use
the equality \thetag{11.2} from this theorem for the further calculations.
Let's transform quadratic in $\bold X$ expression $\bold R(\bold X,\bold Y)
\bold X$ into bilinear one respective to $\bold X$ and $\bold Z$:
$$
\bold R(\bold X,\bold Y)\bold Z+\bold R(\bold Z,\bold Y)\bold X
=\frac{\bold R(\bold X+\bold Z,\bold Y)(\bold X+\bold Z)-\bold R(
\bold X-\bold Z,\bold Y)(\bold X-\bold Z)}{2}.
$$
Applying formula \thetag{11.2} we bring this relationship to the form
$$
\bold R(\bold X,\bold Y)\bold Z+\bold R(\bold Z,\bold Y)\bold X
=\theta(\bold X,\bold Y)\cdot\bold Z+\theta(\bold Z,\bold Y)\cdot
\bold X-2\,\beta(\bold X,\bold Z)\cdot\bold Y.
$$
Here $\bold X$, $\bold Y$, $\bold Z$ are arbitrary vectors. Therefore 
we can get the following relationship for the components of the 
curvature tensor:
$$
\pagebreak
R^k_{sij}+R^k_{isj}=\theta_{ij}\,\delta^k_s+\theta_{sj}\,\delta^k_i-
2\,\beta_{is}\,\delta^k_j.\hskip -3em
\tag11.24
$$
Let's alternate \thetag{11.24} with respect to the pair of indices $i$ 
and $j$:
$$
R^k_{sij}+\frac{R^k_{isj}-R^k_{jsi}}{2}=\tilde\theta_{ij}\,\delta^k_s+
\frac{\theta_{sj}-2\,\beta_{sj}}{2}\,\delta^k_i-\frac{\theta_{si}-2\,
\beta_{si}}{2}\,\delta^k_j.\hskip -3em
\tag11.25
$$
Here $\tilde{\boldsymbol\theta}$ is skew-symmetric part of tensor
$\boldsymbol\theta$. For to transform the left hand side of \thetag{11.25}
we use skew symmetry $R^k_{jsi}=-R^k_{jis}$, $R^k_{sji}=-R^k_{sij}$, and 
the identity $R^k_{sji}+R^k_{isj}+R^k_{jis}=0$, which follows from 
$\Gamma^k_{rs}=\Gamma^k_{sr}$:
$$
\frac{3\,R^k_{sij}}{2}=\tilde\theta_{ij}\,\delta^k_s+\frac{\theta_{sj}
-2\,\beta_{sj}}{2}\,\delta^k_i-\frac{\theta_{si}-2\,\beta_{si}}{2}\,
\delta^k_j.\hskip -3em
\tag11.26
$$
Now let's compare the formula \thetag{10.13} for the case $\dim V\geqslant 
4$ with the formula \thetag{11.26} just derived for the case $\dim V=3$. 
Comparing we see that the structure of components of curvature tensor for 
three-dimensional case is the same as for the cases of higher dimension.
Therefore formula \thetag{10.15} remains true for $n=3$ and we can state
the following theorem similar to theorem~10.2.
\proclaim{Theorem 11.2} Let \thetag{1.1} be a system of equations belonging 
to $m$-th case of intermediate degeneration such that its algebra of point 
symmetries has maximal dimension \thetag{9.1}. Under these assumptions
if $n=3$, then appropriate curvature tensor is expressed through Ricci 
tensor according to the formula \thetag{10.15}.
\endproclaim
\head
12. Structure of curvature tensor for $n=2$.
\endhead
    For $n=2$ due to the inequalities $1\leqslant m\leqslant n-1$
we have only one case of intermediate degeneration with $m=1$.
According to lemma~9.2 tensor $\bold S$ defined in \thetag{4.11}
and skew-symmetric part of Ricci tensor $\tilde\bold R$ in this case
are zero: $\bold S=2\,\tilde\bold R=0$. Curvature tensor is completely 
determined by symmetric part of Ricci tensor $\hat\bold R$:
$$
R^k_{sij}=\hat R_{sj}\,\delta^k_i-\hat R_{si}\,\delta^k_j.\hskip -3em
\tag12.1
$$
In order to prove formula \thetag{12.1} note that for $n=2$ each 
skew-symmetric in pair of indices numeric array is proportional
to skew-symmetric unit matrix
$$
d^{ij}=d_{ij}=\Vmatrix 0 & 1\\ \vspace{1ex}-1 & 0 \endVmatrix.
\hskip -3em
\tag12.2
$$
When applied to curvature tensor, this yields $R^k_{sij}=\rho^k_s\,d_{ij}$.
Hence we can calculate Ricci tensor, which is known to by symmetric by
lemma~9.2:
$$
\hat R_{sj}=R_{sj}=\sum^2_{k=1}\rho^k_s\,d_{kj}.
$$
Matrix \thetag{12.2} is invertible. Therefore we can express $\rho^k_s$
through Ricci tensor:
$$
\pagebreak
\rho^k_s=-\shave{\sum^2_{r=1}}\hat R_{sr}\,d^{rk}.\hskip -3em
\tag12.3
$$
Now we are to substitute \thetag{12.3} into the formula $R^k_{sij}=
\rho^k_s\,d_{ij}$ and use well-known identity $d^{\hskip 0.5pt rk}\,
d_{ij}=\delta^r_i\,\delta^k_j-\delta^k_i\,\delta^r_j$. As a result
we will get \thetag{12.1}.
\head
13. Structure of tensor field $\bold A$.
\endhead
     Let \thetag{1.1} be a system of equations belonging to $m$-th 
case of intermediate degeneration such that its algebra of point 
symmetries has maximal dimension \thetag{9.1}. Choose some fixed
point $p_0$ on $M$ and consider a field of point symmetry $\boldsymbol
\eta$ for \thetag{1.1} vanishing at the point $p_0$. The condition
$L_{\boldsymbol\eta}(\bold A)=0$ from \thetag{2.5} written at the
point $p_0$ for such symmetry field $\boldsymbol\eta$ has the 
following form:
$$
\bold A\,\bold F=\bold F\,\bold A.\hskip -3em
\tag13.1
$$
Here $\bold F=L_{\boldsymbol\eta}$ is linear operator from \thetag{5.7}.
Due to \thetag{9.1} any operator $\bold F$ satisfying \thetag{8.3} 
corresponds to some symmetry field $\boldsymbol\eta\in\Cal G(p_0)$.
In particular, we can choose $\bold F=\bold w\otimes\alpha$, where
$\bold w$ is some arbitrary nonzero vector from the kernel of bilinear 
form $\hat R$ and $\alpha$ is arbitrary linear functional in $V=T_{p_0}(M)$.
Substituting this operator into \thetag{13.1} we obtain the relationship
$$
\alpha(\bold X)\cdot\bold A\bold w=\alpha(\bold A\bold X)\cdot\bold w,
\hskip -3em
\tag13.2
$$
where $\bold X$ is an arbitrary vector of $V$. For $\alpha$ in 
\thetag{13.2} we can choose linear functional such that $\alpha(
\bold X)\neq 0$. Then the equality \thetag{13.2} takes the form
$$
\bold A\bold w=a\cdot\bold w,
$$
where numeric parameter $a$ is defined by formula $a=\alpha(\bold A
\bold X)/\alpha(\bold X)$. Therefore $a$ doesn't depend on the choice
of $\bold w\in\Ker\hat R$. Vector $\bold w\neq 0$ is an eigenvector 
for the operator $\bold A$ corresponding to the eigenvalue $a$. Hence $a$ 
do not depend on the choice of $\bold X$ and $\alpha$ in \thetag{13.2}.
From $\bold A\bold w=a\cdot\bold w$ and \thetag{13.2} we get $\alpha(
\bold A\bold X)=\alpha(a\cdot\bold X)$. Since $\alpha$ is arbitrary
linear functional, we conclude
$$
\bold A\bold X=a\cdot\bold X.\hskip -3em
\tag13.3
$$
The equality \thetag{13.3} fulfilled for arbitrary vector $\bold X\in V$ 
means that $\bold A$ is a scalar operator. It differs from identical 
operator $\id_V$ by a scalar factor: $\bold A=a\cdot\id_V$.\par
    Let's substitute $\bold A=a\cdot\id_V$ into the equation \thetag{2.5},
where point symmetry field $\boldsymbol\eta$ is no longer vanishing at
the point $p_0$:
$$
L_{\boldsymbol\eta}(\bold A)=L_{\boldsymbol\eta}(a)\cdot\id_V=0.
\hskip -3em
\tag13.4
$$
Therefore $L_{\boldsymbol\eta}(a)=0$. From \thetag{9.1} and the estimates
\thetag{8.1} and \thetag{8.9} we get
$$
\dim(\Cal G/\Cal G(p_0))=\dim M=n.
$$
Hence for each point $p_0\in M$ we can find $n$ vector fields from $\Cal G$
whose values at this point are linear independent. Therefore $L_{\boldsymbol
\eta}(a)=0$ \pagebreak implies $a=\const$.
\proclaim{Theorem 13.1} Let \thetag{1.1} be a system of equations 
belonging to $m$-th case of intermediate degeneration such that its algebra 
of point symmetries has maximal dimension \thetag{9.1}. Then $A$ is a 
constant scalar matrix, i\.~e\. $A^i_j=a\,\delta^i_j$ and $a=\const$.
\endproclaim
\head
14. Projectively-euclidean spaces.
\endhead
     For the systems of equations \thetag{1.1} possessing symmetry algebras 
of maximal order \thetag{9.1} matrices $A$ in \thetag{1.1} are the same in 
all cases of intermediate degeneration: $A^i_j=a\,\delta^i_j$. Therefore 
structure of such systems of equations are completely defined by components 
of affine connection $\Gamma^j_{rs}$ in \thetag{1.1}.\par
     Let $\bold R$ be curvature tensor corresponding to the connection 
$\Gamma$ in \thetag{1.1}. According to theorems~10.2 and 11.2 for $n
\geqslant 3$ its components are given by formula \thetag{10.15}. In order 
to simplify notations let's use formula \thetag{10.13} instead of 
\thetag{10.15}. Tensor fields $\boldsymbol\beta$ and $\boldsymbol\sigma$ 
in \thetag{10.13} are expressed through Ricci tensor according to 
\thetag{10.14}. Taking into account the relationship $\boldsymbol\sigma=
-2\,\tilde{\boldsymbol\beta}$, which follows from \thetag{10.14}, we can
bring \thetag{10.13} to the following rather simple form:
$$
R^k_{sij}=(\beta_{ji}-\beta_{ij})\,\delta^k_s+\beta_{js}\,\delta^k_i
-\beta_{is}\,\delta^k_j.\hskip -3em
\tag14.1
$$
Let's substitute \thetag{14.1} into the well-known Bianchi-Padov
identity (see \cite{17}, \cite{22}):
$$
\nabla_iR^k_{sjq}+\nabla_jR^k_{sqi}+\nabla_qR^k_{sij}=0.
$$
If $n\geqslant 3$, this substitution results in the equation 
for tensor field $\boldsymbol\beta$:
$$
\nabla_i\beta_{js}=\nabla_j\beta_{is}.\hskip -3em
\tag14.2
$$
According to \cite{22} (see chapter 4, \S\,47) we write the following
differential equation for some covector field $\bold u$:
$$
\nabla_iu_j=\beta_{ij}+u_i\,u_j,\quad i,j=1,\,\ldots,\,n.
\hskip -3em
\tag14.3
$$
If we express covariant derivatives through partial derivatives, we see
that \thetag{14.3} is a complete system of Pfaff equations:
$$
\frac{\partial u_j}{\partial y^i}=\beta_{ij}+u_i\,u_j+\sum^n_{k=1}
\Gamma^k_{ij}\,u_k.\hskip -3em
\tag14.4
$$
By direct calculations one can find that the relationships \thetag{14.1} 
and \thetag{14.2} provide complete compatibility for Pfaff equations
\thetag{14.4}. This in turn means that for differential equations 
\thetag{14.3} the Cauchy problem 
$$
u_j\,\hbox{\vrule height 6pt depth 8pt width 0.5pt}_{\,p=p_0}
=u_j(0)\hskip -3em
\tag14.5
$$
is solvable in some neighborhood of the point $p=p_0$ in $M$ for arbitrary
initial values $u_j(0)$ in \thetag{14.5}.\par
     Let $\bold u$ be a covector field obtained by solving Cauchy problem
\thetag{14.5} for the equations \thetag{14.3}. We shall use this field to
construct another connection $\bar\Gamma$:
$$
\bar\Gamma^k_{rs}=\Gamma^k_{rs}+u_r\,\delta^k_s+u_s\,\delta^k_r.\hskip -3em
\tag14.6
$$
Let's calculate curvature tensor of new connection $\bar\Gamma$ on the
base of formula \thetag{2.8}. Taking in to account \thetag{14.1} and 
\thetag{14.3} we get $\bar R^k_{sij}=0$. This means that \thetag{14.1}
is flat (euclidean) affine connection.\par
     Transformation of connection components $\bar\Gamma^k_{rs}\to 
\Gamma^k_{rs}=\bar\Gamma^k_{rs}-u_r\,\delta^k_s-u_s\,\delta^k_r$ is
called {\it projective transformation} (see \cite{22}). Affine connection 
$\Gamma$ obtained from euclidean connection $\bar\Gamma$ by this 
transformation is called {\it projectively-euclidean connection}. Manifolds 
equipped with projectively-euclidean connection are called 
{\it projecti\-vely-euclidean spaces}. For euclidean connection \thetag{14.6} 
one can find local coordinates $y^1,\,\ldots,\,y^n$ such that $\bar
\Gamma^k_{rs}=0$. Relative to projectively-euclidean connection $\Gamma$ 
these local coordinates are called {\it projectively-euclidean coordinates}. 
In  projectively-euclidean coordinates for $\Gamma^k_{rs}$ we have
$$
\Gamma^k_{rs}=-(u_r\,\delta^k_s+u_s\,\delta^k_r).\hskip -3em
\tag14.7
$$
\proclaim{Theorem 14.1} Let \thetag{1.1} be a system of equations belonging 
to $m$-th case of intermediate degeneration such that its algebra of point 
symmetries has maximal dimension \thetag{9.1}. Under these assumptions
if $n\geqslant 3$, then there is a point transformation \thetag{1.4}
bringing these equations to the form
$$
\frac{\partial y^i}{\partial\,t}=a\,\frac{\partial^2 y^i}{\partial x^2}
-2\,a\left(\,\shave{\sum^n_{r=1}}u_r\,\frac{\partial y^r}{\partial x}
\right)\frac{\partial y^i}{\partial x},\quad i=1,\,\ldots,\,n,
$$
where $a=\const$ and $u_r=u_r(y^1,\ldots,y^n)$ are components of some
covector field $\bold u$.
\endproclaim
\head
15. Structure of tensor $\bold Q=\nabla\bold R$.
\endhead
     Denote by $\bold Q$ covariant differential of Ricci tensor: $\bold Q
=\nabla\bold R$ (in local coordinates $Q_{ijk}=\nabla_iR_{jk}$). $\bold Q$ 
is a tensor field of type $(0,3)$ belonging to the algebra $\Cal R$ from 
\thetag{4.10}. According to the theorem~4.3 we have
$$
L_{\boldsymbol\eta}(\bold Q)=0\hskip -3em
\tag15.1
$$
for any point symmetry field $\boldsymbol\eta$ of the system of 
equations \thetag{1.1}. Let's take and fix some point $p=p_0$ on 
$M$. Tensor $\bold Q$ defines trilinear form on $V=T_{p_0}(M)$:
$$
Q(\bold X,\bold Y,\bold Z)=\sum^n_{i=1}\sum^n_{j=1}\sum^n_{k=1}
\nabla_iR_{jk}\,Z^i\,X^j\,Y^k.\hskip -3em
\tag15.2
$$
From \thetag{15.1} we get the following equality for trilinear 
form \thetag{15.2}:
$$
Q(\bold F\bold X,\bold Y,\bold Z)+Q(\bold X,\bold F\bold Y,\bold Z)
+Q(\bold X,\bold Y,\bold F\bold Z)=0.\hskip -3em
\tag15.3
$$
Here $\bold F$ is linear operator \thetag{5.7} defined by point symmetry
field $\boldsymbol\eta$ vanishing at $p_0$.
\proclaim{Lemma 15.1} Let \thetag{1.1} be a system of equations belonging 
to $m$-th case of intermediate degeneration such that its algebra of point 
symmetries has maximal dimension \thetag{9.1}. Under these assumptions
if $n\geqslant 3$, then 
$$
Q(\bold w,\bold Y,\bold Z)=Q(\bold X,\bold w,\bold Z)=Q(\bold X,\bold Y,
\bold w)=0,\hskip -3em
\tag15.4
$$
where $\bold w$ is a vector from kernel of quadratic form $\hat R$ given
by symmetric part of Ricci tensor and $\bold X$, $\bold Y$, $\bold Z$ are
arbitrary three vectors in $V$.
\endproclaim
\demo{Proof} Relying on the results of section 8, we consider linear
operator $\bold F=\bold w\otimes\alpha$, where $\alpha$ is an arbitrary
linear functional in $V$. Substituting $\bold F$ into \thetag{15.3} we get
$$
\alpha(\bold X)\,Q(\bold w,\bold Y,\bold Z)+\alpha(\bold Y)\,Q(\bold X,
\bold w,\bold Z)+\alpha(\bold Z)\,Q(\bold X,\bold Y,\bold w)=0.\hskip -3em
\tag15.5
$$
In the space $V$ of dimension $n\geqslant 3$ for any two vectors $\bold Y$ 
and $\bold Z$ one can find third vector $\bold X$ and a linear functional
$\alpha$ such that 
$$
\xalignat 3
&\alpha(\bold X)=1,&&\alpha(\bold Y)=0,&&\alpha(\bold Z)=0.
\endxalignat
$$
Then from \thetag{15.5} we get $Q(\bold w,\bold Y,\bold Z)=0$, thus
proving first relationship \thetag{15.4}. Other two relationships
are proved in a similar way.\qed\enddemo
\proclaim{Lemma 15.2} Proposition of lemma~15.1 remains true for the 
dimension $n=2$.
\endproclaim
\demo{Proof} For $n=2$ from inequalities $1\leqslant m\leqslant n-1$
we get $m=1$. Due to lemma~9.2 skew-symmetric part of Ricci tensor is 
zero: $\tilde\bold R=0$. Therefore tensor $\bold Q$ is symmetric 
respective to last pair of indices. For trilinear form \thetag{15.2} 
we obtain
$$
Q(\bold X,\bold Y,\bold Z)=Q(\bold Y,\bold X,\bold Z).\hskip -3em
\tag15.6
$$
Define a form $\phi(\bold X,\bold Z)=Q(\bold X,\bold X,\bold Z)$, which
is linear respective to $\bold Z$ and quadratic respective to $\bold X$.
For this form from \thetag{15.5} we derive
$$
2\,\alpha(\bold X)\,Q(\bold w,\bold X,\bold Z)
+\alpha(\bold Z)\,\phi(\bold X,\bold w)=0.\hskip -3em
\tag15.7
$$
For any vector $\bold X$ in two-dimensional space one can find another
vector $\bold Z$ and linear functional $\alpha$ such that $\alpha(\bold X)
=0$ and $\alpha(\bold Z)=1$. Due to \thetag{15.7} this yields $\phi(\bold X,
\bold w)=0$. Form $Q$ can be recovered by form $\phi$ due to its symmetry
(see formula \thetag{15.6} above). We have an obvious relationship
$$
Q(\bold X,\bold Y,\bold Z)=\frac{\phi(\bold X+\bold Y,\bold Z)
-\phi(\bold X-\bold Y,\bold Z)}{4}.\hskip -3em
\tag15.8
$$
Substituting $\bold Z=w$ into the formula \thetag{15.8} and taking into
account $\phi(\bold X,\bold w)=0$, we get one of the required relationships
$Q(\bold X,\bold Y,\bold w)=0$.\par
     In order to prove rest two relationships in \thetag{15.4} we take into
account that $Q(\bold X,\bold Y,\bold w)=0$. This brings \thetag{15.5} to
the following form:
$$
\pagebreak
\alpha(\bold X)\,Q(\bold w,\bold Y,\bold Z)+\alpha(\bold Y)\,Q(\bold X,
\bold w,\bold Z)=0.\hskip -3em
\tag15.9
$$
For any vector $\bold X$ in two-dimensional space one can find another
vector $\bold Y$ and linear functional $\alpha$ such that $\alpha(\bold X)
=0$ and $\alpha(\bold Y)=1$. Therefore from \thetag{15.9} we get $Q(\bold 
X,\bold w,\bold Z)=0$. The last relationship $Q(\bold w,\bold Y,\bold Z)=0$
follows from previous one due to \thetag{15.6}. Lemma~15.2 is proved.
\qed\enddemo
\proclaim{Lemma 15.3} Let \thetag{1.1} be a system of equations belonging 
to $m$-th case of intermediate degeneration such that its algebra of point 
symmetries has maximal dimension \thetag{9.1}. Under these assumptions
if $n\geqslant 3$, then $\bold Q=\nabla\bold R=0$.
\endproclaim
\demo{Proof} Due to lemmas~15.1 and 15.2 trilinear form $Q$ from
\thetag{15.2} is zero in the kernel $W=\Ker\hat R$ of bilinear form
$\hat R$ from \thetag{8.3}. This circumstance makes possible to defines 
an induced trilinear form $Q$ in factorspace $V/W$ where induces form
$\hat R$ is nondegenerate. We define this form by the relationship
$$
Q(\hat\bold X,\hat\bold Y,\hat\bold Z)=Q(\bold X,\bold Y,
\bold Z),\hskip -3em
\tag15.10
$$
where $\hat\bold X=\Cl_{\sssize W}(\bold X)$, $\hat\bold Y=\Cl_{\sssize W}
(\bold Y)$, and $\hat\bold Z=\Cl_{\sssize W}(\bold Z)$ are cosets of three
vectors $\bold X$, $\bold Y$, and $\bold Z$ respective to subspace $W$.
According to lemma~9.2 for $m\geqslant 3$ skew-symmetric part of Ricci
tensor is zero. Therefore the relationship \thetag{15.6} holds and 
hence form \thetag{15.10} in factorspace is symmetric:
$$
Q(\hat\bold X,\hat\bold Y,\hat\bold Z)=Q(\hat\bold Y,\hat\bold X,
\hat\bold Z).\hskip -3em
\tag15.11
$$
By applying the complexification procedure to the factorspace $V/W$, 
if necessary, the matrix of nondegenerate symmetric bilinear form 
$\hat R$ from \thetag{8.7} can be brought to the unit matrix at the 
expense of proper choice of base $\bold e_1,\,\ldots,\,\bold e_m$ 
in factorspace $V/W$ (see relationships \thetag{9.7}). Now consider
an operator $\hat\bold F\in\End(V/W)$ defined by \thetag{9.8}. For
this operator we have the relationship $\hat R(\hat\bold F\hat\bold X,
\hat\bold Y)=-\hat R(\hat\bold X,\hat\bold F\hat\bold Y)$, which
follows from \thetag{9.7} and \thetag{9.8}. Due to \thetag{9.1}
operator $\hat\bold F$ should satisfy the equality
$$
Q(\hat\bold F\hat\bold X,\hat\bold Y,\hat\bold Z)+Q(\hat\bold X,
\hat\bold F\hat\bold Y,\hat\bold Z)+Q(\hat\bold X,\hat\bold Y,
\hat\bold F\hat\bold Z)=0,\hskip -3em
\tag15.12
$$
which is obtained from \thetag{15.3} by factorization respective
to subspace $W=\Ker\hat R$.\par
    Let's substitute $\bold e_k$, $\bold e_q$, and $\bold e_i$ for
$\hat\bold X$, $\hat\bold Y$, and $\hat\bold Z$ into \thetag{15.12} 
in various combinations. Dimension $\dim(V/W)=m\geqslant 3$ is enough 
to choose indices $i$, $k$, and $q$ to be mutually distinct. If we 
take into account \thetag{9.8} and substitute $\hat\bold X=\bold e_k$, 
$\hat\bold Y=\bold e_k$, $\hat\bold Z=\bold e_i$ into \thetag{15.12},
we get $Q(\bold e_q,\bold e_k,\bold e_i)=-Q(\bold e_k,\bold e_q,\bold 
e_i)$. From the formula \thetag{15.11} in turn we get $Q(\bold e_q,
\bold e_k,\bold e_i)=Q(\bold e_k,\bold e_q,\bold e_i)$. Therefore
$$
Q_{qki}=Q(\bold e_q,\bold e_k,\bold e_i)=0\text{, \ for \ }
q\neq k\neq i\neq q.\hskip -3em
\tag15.13
$$
Now let's make other three substitutions into the formula
\thetag{15.12}. Their results are expressed by the following
implications:
$$
\align
&\hat\bold X=\bold e_i,\ \hat\bold Y=\bold e_i,\ \hat\bold Z
=\bold e_k\text{\ \ implies \ }Q(\bold e_i,\bold e_i,\bold e_q)=0,\\
&\hat\bold X=\bold e_i,\ \hat\bold Y=\bold e_k,\ \hat\bold Z
=\bold e_i\text{\ \ implies \ }Q(\bold e_i,\bold e_q,\bold e_i)=0,\\
\displaybreak
&\hat\bold X=\bold e_k,\ \hat\bold Y=\bold e_i,\ \hat\bold Z
=\bold e_i\text{\ \ implies \ }Q(\bold e_q,\bold e_i,\bold e_i)=0.
\endalign
$$
Let's combine all these results into one formula
$$
Q(\bold e_i,\bold e_i,\bold e_q)=Q(\bold e_i,\bold e_q,\bold e_i)=
Q(\bold e_q,\bold e_i,\bold e_i)=0\text{\ \ for \ }i\neq q.\hskip -3em
\tag15.14
$$
And finally, let's substitute $\hat\bold X=\bold e_k$, $\hat\bold Y=
\bold e_q$, $\hat\bold Z=\bold e_q$ into \thetag{15.12}. Here we
get $Q(\bold e_q,\bold e_q,\bold e_q)-Q(\bold e_k,\bold e_k,\bold e_q)
-Q(\bold e_k,\bold e_q,\bold e_k)=0$. If we take into account 
\thetag{15.14}, we can bring this equality to the following form
$$
Q(\bold e_q,\bold e_q,\bold e_q)=0.\hskip -3em
\tag15.15
$$
Now summarizing \thetag{15.13}, \thetag{15.14}, and \thetag{15.15} we
see that all components of trilinear form \thetag{15.10} in the base
$\bold e_1,\,\ldots,\,\bold e_m$ of factorspace $V/W$ are zero. This
implies that $Q$ iz zero in $V/W$, and hence, due to \thetag{15.10}
$Q$ iz zero in $V$. Therefore $\nabla\bold R=\bold Q=0$. Lemma is
proved.\qed\enddemo
\proclaim{Lemma 15.4} Proposition of lemma~15.3 remains true for the 
case $m=2$.
\endproclaim
\demo{Proof} For $m=2$ we are restricted by two non-coinciding values
of indices $q\neq k$ for numerating base vectors in $V/W$. As above
we suppose the relationships \thetag{9.7} to be fulfilled for the
base $\bold e_1,\,\bold e_2$ in $V/W$ and consider an operator $\hat
\bold F\in\End(V,W)$ defined by it action on the base vectors:
$$
\hat\bold F(\bold e_i)=\cases\hphantom{-}\bold e_q &\text{for \ }
i=k\neq q,\\
-\bold e_k &\text{for \ }i=q\neq k.\endcases\hskip -3em
\tag15.16
$$
Relying on \thetag{15.16} we make the following substitutions into 
the formula \thetag{15.12}:
$$
\align
&\hat\bold X=\bold e_q,\ \hat\bold Y=\bold e_k,\ \hat\bold Z
=\bold e_q\text{\ \ implies \ }Q(\bold e_q,\bold e_q,\bold e_q)=
Q(\bold e_k,\bold e_k,\bold e_q)+Q(\bold e_q,\bold e_k,\bold e_k),\\
\vspace{1ex}
&\hat\bold X=\bold e_k,\ \hat\bold Y=\bold e_q,\ \hat\bold Z
=\bold e_q\text{\ \ implies \ }Q(\bold e_q,\bold e_q,\bold e_q)=
Q(\bold e_k,\bold e_k,\bold e_q)+Q(\bold e_k,\bold e_q,\bold e_k),\\
\vspace{1ex}
&\hat\bold X=\bold e_q,\ \hat\bold Y=\bold e_q,\ \hat\bold Z
=\bold e_k\text{\ \ implies \ }Q(\bold e_q,\bold e_q,\bold e_q)=
Q(\bold e_k,\bold e_q,\bold e_k)+Q(\bold e_q,\bold e_k,\bold e_k).
\endalign
$$
Let's subtract first of the above relationships from the second one
and from the third one. The result can be written as
$$
Q(\bold e_q,\bold e_k,\bold e_k)=Q(\bold e_k,\bold e_q,\bold e_k)
=Q(\bold e_k,\bold e_k,\bold e_q).\hskip -3em
\tag15.17
$$
To the contrary, if we add all of them and take into account \thetag{15.17},
we get
$$
Q(\bold e_q,\bold e_q,\bold e_q)=2\,Q(\bold e_k,\bold e_q,\bold e_k).
\hskip -3em
\tag15.18
$$
Now let's substitute $\hat\bold X=\bold e_k$, $\hat\bold Y=\bold e_k$,
and $\hat\bold Z=\bold e_k$ into \thetag{15.12}. This yields
$$
Q(\bold e_q,\bold e_k,\bold e_k)+Q(\bold e_k,\bold e_q,\bold e_k)
+Q(\bold e_k,\bold e_k,\bold e_q)=0.\hskip -3em
\tag15.19
$$
From \thetag{15.17}, \thetag{15.18}, and \thetag{15.19} we see that
all components of trilinear form $Q$ in the base $\bold e_1,\,\bold e_2$ 
of factorspace $V/W$ are zero:
$$
Q(\bold e_q,\bold e_k,\bold e_k)=Q(\bold e_k,\bold e_q,\bold e_k)
=Q(\bold e_k,\bold e_k,\bold e_q)=Q(\bold e_q,\bold e_q,\bold e_q)=0
\text{\ \ for \ }q\neq k.
$$
Then from \thetag{15.10} we get $Q=0$, thus completing the the proof of
lemma~15.4.\qed\enddemo
     For $m=1$ tensor $\bold Q=\nabla\bold R$ shouldn't vanish. However,
in this case factorspace $V/W$ is one-dimensional. Trilinear form
\thetag{15.10} in $V/W$ has only one component $Q_{111}=Q(\bold e_1,
\bold e_1,\bold e_1)$. Therefore trilinear form $Q$ in the initial space
$V$ is determined by some linear form $u$ up to a scalar factor:
$$
Q(\bold X,\bold Y,\bold Z)\sim u(\bold X)\,u(\bold Y)\,u(\bold Z).
\hskip -3em
\tag15.20
$$
For tensor field $\bold Q=\nabla\bold R$ formula \thetag{15.20} is 
written as $\bold Q=\nabla\bold R\sim \bold u\otimes\bold u\otimes
\bold u$. This means that trilinear form \thetag{15.20} is completely
symmetric. For covariant differential of Ricci tensor this yields
the following symmetry:
$$
\nabla_i R_{jk}=\nabla_j R_{ik}.\hskip -3em
\tag15.21
$$
Ricci tensor itself for $m=1$ is also symmetric: $\bold R=\hat\bold R$
and $S=2\,\tilde\bold R=0$, as stated in lemma~9.2. For bilinear
form $\hat R$ we have the relationship
$$
\hat R(\bold X,\bold Y)\sim u(\bold X)\,u(\bold Y),\hskip -3em
\tag15.22
$$
similar to \thetag{15.20}. Due to \thetag{15.22} linear form
$u(\bold X)$ should vanish when applied to the vector $\bold X$
in the kernel $W=\Ker\hat R$.
\head
16. Two-dimensional projectively-euclidean spaces.
\endhead
     For $n=2$ we have only one case of intermediate degeneration
with $m=1$. Suppose that system of equations \thetag{1.1} belongs
to this case and suppose that it has an algebra of point symmetries
of maximal dimension \thetag{9.1}. Ricci tensor for the connection
$\Gamma$ defined by this system of equations is symmetric ($\bold R
=\hat\bold R$, $\tilde\bold R=0$). Curvature tensor has the form 
\thetag{12.1}, which is same as \thetag{10.15} if we take into
account $n=2$ and $\tilde\bold R=0$. Therefore we can transform
\thetag{12.1} to \thetag{14.1} if we denote $\beta_{ij}=\hat R_{ij}
=R_{ij}$. The equations \thetag{14.2} for the field $\boldsymbol\beta$ 
in two-dimensional case are not immediate consequence of \thetag{12.1}
(Bianchi-Padov identity for $n=2$ is trivial). However, they follow from 
\thetag{9.1} due to lemma~15.2 which leads to \thetag{15.21}. Since the 
equations \thetag{14.2} are fulfilled, we can repeat arguments from 
section~14 for two-dimensional case $n=2$. This results in the following 
theorem similar to theorem~14.1.
\proclaim{Theorem 16.1} Let \thetag{1.1} be a system of equations belonging 
to $m$-th case of intermediate degeneration \pagebreak such that its 
algebra of point symmetries has maximal dimension \thetag{9.1}. Under these 
assumptions if $n=2$, then there is a point transformation \thetag{1.4} 
bringing these equations to the form
$$
\frac{\partial y^i}{\partial\,t}=a\,\frac{\partial^2 y^i}{\partial x^2}
-2\,a\left(\,\shave{\sum^2_{r=1}}u_r\,\frac{\partial y^r}{\partial x}
\right)\frac{\partial y^i}{\partial x},\quad i=1,\,2,
$$
where $a=\const$ and $u_r=u_r(y^1,y^2)$ are components of some
covector field $\bold u$.
\endproclaim
\head
17. Separation of variables.
\endhead
     Due to results of section~16 we can make further consideration
common for all cases of intermediate degeneration $1\leqslant m\leqslant 
n-1$ and for all dimensions $n\geqslant 2$. In order to construct
projectively-euclidean coordinates in section~14 we chose some solution
of the system of Pfaff equations \thetag{14.3} defining some covector
field $\bold u$ on the manifold $M$. Now we shall make the choice of
field $\bold u$ more specific. In order to do it we complete the
equations \thetag{14.3} with the equations for some vector field
$\boldsymbol\xi$:
$$
\cases
\nabla_iu_k=\beta_{ik}+u_i\,u_k,\\
\vspace{1ex}
\nabla_i\xi^k=-u_i\,\xi^k.
\endcases\hskip -3em
\tag17.1
$$
Here $i,k=1,\,\ldots,\,n$. The equations \thetag{17.1} form complete
system of Pfaff equations. However, in contrast to \thetag{14.3} this
system of Pfaff equations is not completely compatible. In order to
make it compatible we impose the following restrictions for the
variables $u_1,\,\ldots,\,u_n$ and $\xi^1,\,\ldots,\,\xi^n$:
$$
\xalignat 2
&\sum^n_{r=1}\hat\beta_{ir}\,\xi^r=0,
&&\sum^n_{r=1}u_r\,\xi^r=0\hskip -3em
\tag17.2
\endxalignat
$$
(brief introduction into the theory of Pfaff equations with restrictions
can be found in \cite{23}, appendix~B). Let's check the compatibility of
\thetag{17.1} due to restrictions \thetag{17.2}. For to do it we calculate 
the commutator of covariant derivatives according to the equations 
\thetag{17.1}. This is implemented by the following two equations:
$$
\aligned
-&\sum^n_{s=1}R^s_{kij}\,u_s=[\nabla_i,\nabla_j]u_k=\nabla_i(\beta_{jk}
+u_j\,u_k)-\nabla_j(\beta_{ik}+u_i\,u_k),\\
&\sum^n_{s=1}R^k_{sij}\,\xi^s=[\nabla_i,\nabla_j]\xi^k=\nabla_i(-u_j\,
\xi^k)-\nabla_j(u_i\,\xi^k).
\endaligned\hskip -3em
\tag17.3
$$
First equation \thetag{17.3} is an identity due to \thetag{14.1} and
\thetag{14.2}, here we need no restrictions. Second one is reduced to
$$
\sum^n_{s=1}\beta_{js}\,\xi^s\,\delta^k_i-\sum^n_{s=1}\beta_{is}\,\xi^s
\,\delta^k_j=0.\hskip -3em
\tag17.4
$$
The equality \thetag{17.4} holds due to the first equation of
restrictions \thetag{17.2}. Indeed, \thetag{17.2} means that
vector $\boldsymbol\xi$ belongs to the kernel of quadratic form
$\hat R$. By lemma~9.1 it follows that the kernel of skew-symmetric
part of the form $\beta$ contains the kernel of its symmetric part:
$\Ker\tilde\beta=\Ker S=\Ker\tilde R\supset\Ker\hat R$.\par
     Now let's prove that the restrictions \thetag{17.2}  are  in 
agreement
with the equations \thetag{17.1}. In order to do it let's calculate
covariant derivatives of left hand sides of \thetag{17.2} and equate 
them to zero:
$$
\align
0&=\sum^n_{r=1}\nabla_i(\hat\beta_{jr}\,\xi^r)=\sum^n_{r=1}\nabla_i
\hat\beta_{jr}\,\xi^r-\sum^n_{r=1}u_i\,\hat\beta_{jr}\,\xi^r,
\hskip -3em
\tag17.5\\
0&=\sum^n_{r=1}\nabla_i(u_r\,\xi^r)=\sum^n_{r=1}\beta_{ir}\,\xi^r.
\hskip -3em
\tag17.6
\endalign
$$
The equality \thetag{17.6} is fulfilled due to the first equation
of restrictions \thetag{17.2}. The equality \thetag{17.5} also
is provided by this equation and lemmas~15.1 and 15.2.\par 
     Compatibility of the Pfaff equations \thetag{17.1} restricted
by the equations \thetag{17.2} means that for these equations the 
Cauchy problem
$$
\xalignat 2
&u_k\,\hbox{\vrule height 6pt depth 8pt width 0.5pt}_{\,p=p_0}
=u_k(0),
&&\xi^k\,\hbox{\vrule height 6pt depth 8pt width 0.5pt}_{\,p=p_0}
=\xi^k(0)\hskip -3em
\tag17.7
\endxalignat
$$
is uniquely solvable in some neighborhood of the point $p_0$
for arbitrary initial values satisfying the restrictions
\thetag{17.2}. Solution of this Cauchy problem satisfies the
equations of restrictions \thetag{17.2} in whole neighborhood 
where it is defined.\par
    Let's choose one set of initial values for the components of
covector $\bold u$ and $n-m$ sets of initial values  for the
components of vector $\boldsymbol\xi$. Vectors $\boldsymbol\xi_1(0),
\,\ldots,\,\boldsymbol\xi_{n-m}(0)$ should be chosen so that they
are linearly independent and all belong to the kernel of the form
$\hat R$. This provides first condition \thetag{17.2}. Covector
$\bold u(0)$ should be chosen such that it vanishes being contracted
with the vectors of the subspace $\Ker\hat R$. This provides second
condition \thetag{17.2}. Pfaff equation \thetag{17.1} respective to
$\bold u$ do not depend on components of $\boldsymbol\xi$. Therefore
by solving $n-m$ Cauchy problems we obtain $n-m$ vector fields
$\boldsymbol\xi_1,\,\ldots,\,\boldsymbol\xi_{n-m}$ and only one
covector field $\bold u$. Vector fields $\boldsymbol\xi_1,\,\ldots,\,
\boldsymbol\xi_{n-m}$ constitute moving frame forming the base of
subspace $\Ker\hat R$ at each point on $M$. Covector field $\bold u$
vanishes being contracted with vectors of frame $\boldsymbol\xi_1,\,
\ldots,\,\boldsymbol\xi_{n-m}$. We express this circumstance as 
$$
u(\boldsymbol\xi_1)=\ldots=u(\boldsymbol\xi_{n-m})=0.\hskip -3em
\tag17.8
$$
Using \thetag{17.1} and the relationships \thetag{17.8} we can calculate 
covariant derivatives of the fields $\boldsymbol\xi_1,\,\ldots,\,
\boldsymbol\xi_{n-m}$, and $\bold u$ along the vectors of the frame 
$\boldsymbol\xi_1,\,\ldots,\,\boldsymbol\xi_{n-m}$:
$$
\xalignat 2
&\nabla_{\boldsymbol\xi_i}\boldsymbol\xi_j=0,
&&\nabla_{\boldsymbol\xi_i}\bold u=0.\hskip -3em
\tag17.9
\endxalignat
$$
Then from \thetag{17.9} we derive the commutators
$$
[\boldsymbol\xi_i,\,\boldsymbol\xi_j]=\nabla_{\boldsymbol\xi_i}
\boldsymbol\xi_{j}-\nabla_{\boldsymbol\xi_j}\boldsymbol\xi_{i}=0.
\hskip -3em
\tag17.10
$$
Following theorem is an immediate consequence of the relationships
\thetag{17.10}.
\proclaim{Theorem 17.1} Let \thetag{1.1} be a system of equations belonging 
to $m$-th case of intermediate degeneration such that its algebra of point 
symmetries has maximal dimension \thetag{9.1}. Under these assumptions 
kernels of bilinear forms $\hat R$ defined by symmetric part of Ricci tensor 
$\hat\bold R$ form involutive $(n-m)$-dimensional distribution on $M$,
which is integrable by Frobenius theorem.
\endproclaim
     Retaining field $\bold u$ defined by \thetag{17.1} and \thetag{17.7} 
as described above let's consider another system of Pfaff equations for
the components of vector field $\bold X$:
$$
\nabla_iX^k=-u_i\,X^k-\sum^n_{r=1}u_r\,X^r\,\delta^k_i,\quad i,
k=1,\,\ldots,\,n.\hskip -3em
\tag17.11
$$
Pfaff equations \thetag{17.11} are completely compatible without 
imposing any restrictions. For the equations \thetag{17.11} the 
Cauchy problem similar to \thetag{17.7} is solvable for arbitrary 
initial values in it:
$$
X^k\,\hbox{\vrule height 6pt depth 8pt width 0.5pt}_{\,p=p_0}
=X^k(0).\hskip -3em
\tag17.12
$$
Let's choose $m$ sets of initial values in \thetag{17.12} such that
cosets of corresponding vectors $\bold X_1(0),\,\ldots,\,\bold X_m(0)$
respective to subspace $W=\Ker\hat R$ are linear independent in
factorspace $V/W$. By solving Cauchy problems \thetag{17.12} for 
\thetag{17.11} with these sets of initial values we obtain $m$ vector
fields $\bold X_1,\,\ldots,\,\bold X_m$ that in aggregate with 
$\boldsymbol\xi_1,\,\ldots,\,\boldsymbol\xi_{n-m}$ form $n$-dimensional
moving frame on $M$. Denote by $u(\bold X_1),\,\ldots,\,u(\bold X_m)$
the results of contraction of covector field $\bold u$ with vector
fields $\bold X_1,\,\ldots,\,\bold X_m$. Then, by direct calculations,
from \thetag{17.11} we get
$$
\nabla_{\bold X_i}\bold X_j=-u(\bold X_i)\bold X_j-u(\bold X_j)
\bold X_i,\hskip -3em
\tag17.13
$$
Analogously from \thetag{17.1} and \thetag{17.11} we derive
$$
\xalignat 2
&\nabla_{\bold X_i}\boldsymbol\xi_j=-u(\bold X_i)\boldsymbol\xi_j,
&&\nabla_{\boldsymbol\xi_i}\bold X_j=-u(\bold X_j)\boldsymbol\xi_i.
\hskip -3em
\tag17.14
\endxalignat
$$
By \thetag{17.14} we find cross commutators of $\bold X_1,\,\ldots,\,
\bold X_m$ and $\boldsymbol\xi_1,\,\ldots,\,\boldsymbol\xi_{n-m}$:
$$
[\boldsymbol\xi_i,\,\bold X_j]=\nabla_{\boldsymbol\xi_i}\bold X_j-
\nabla_{\bold X_j}\boldsymbol\xi_i=0.\hskip -3em
\tag17.15
$$
Moreover, vector fields $\bold X_1,\,\ldots,\,\bold X_m$ are commutating
with each other:
$$
[\bold X_i,\,\bold X_j]=\nabla_{\bold X_i}\bold X_j-
\nabla_{\bold X_j}\bold X_i=0.\hskip -3em
\tag17.16
$$
This follows from \thetag{17.13}. Further we shall calculate derivatives 
of the functions $u(\bold X_1),\,\ldots,\,u(\bold X_m)$ along vector 
fields $\boldsymbol\xi_1,\,\ldots,\,\boldsymbol\xi_{n-m}$:
$$
\pagebreak 
\nabla_{\boldsymbol\xi_i}(u(\bold X_j))=C(\nabla_{\boldsymbol\xi_i}
\bold u\otimes\bold X_j+\bold u\otimes\nabla_{\boldsymbol\xi_i}
\bold X_j)=0.\hskip -3em
\tag17.17
$$\par
     From \thetag{17.10}, \thetag{17.15}, and \thetag{17.16} we conclude
that moving frame $\bold X_1,\,\ldots,\,\bold X_m$, $\boldsymbol\xi_1,\,
\ldots,\,\boldsymbol\xi_{n-m}$ consists of $n$ mutually commutating
vector fields. Therefore one can find local coordinates $y^1,\,\ldots,\,
y^n$ whose frame of coordinate tangent fields coincides with $\bold X_1,
\,\ldots,\,\bold X_m,\boldsymbol\xi_1,\,\ldots,\,\boldsymbol\xi_{n-m}$. 
Formulas \thetag{17.9}, \thetag{17.13}, and \thetag{17.14} 
determine connection components in these local coordinates:
$$
\Gamma^k_{rs}=-(u_r\,\delta^k_s+u_s\,\delta^k_r).\hskip -3em
\tag17.18
$$
Here $u_k=0$ for $k>m$ and $u_k=u_k(y^1,\ldots,y^m)$ for $1\leqslant k
\leqslant m$. This follows from \thetag{17.8} and \thetag{17.17}.
Comparing \thetag{17.18} and \thetag{14.7} we state the theorem
strengthening the result of theorem~14.1. 
\proclaim{Theorem 17.2} Let \thetag{1.1} be a system of equations belonging 
to $m$-th case of intermediate degeneration such that its algebra of point 
symmetries has maximal dimension \thetag{9.1}. Under these assumptions 
there is a point transformation \thetag{1.4} bringing these equations to the 
form
$$
\frac{\partial y^i}{\partial\,t}=a\,\frac{\partial^2 y^i}{\partial x^2}
-2\,a\left(\,\shave{\sum^m_{r=1}}u_r\,\frac{\partial y^r}{\partial x}
\right)\frac{\partial y^i}{\partial x},\quad i=1,\,\ldots,\,n,\hskip -3em
\tag17.19
$$
where $a=\const$ and $u_r=u_r(y^1,\ldots,y^m)$ are components of some
covector field $\bold u$.
\endproclaim
     Let $m=2$. From lemmas~15.3 and 15.4 we know that $\nabla_q\beta_{ij}
=0$. Therefore $\nabla_k\nabla_q\beta_{ij}=0$. Let's write the identity
$$
-\sum^n_{s=1}R^s_{ikq}\,\beta_{sj}-\sum^n_{s=1}R^s_{jkq}\,\beta_{is}
=[\nabla_k,\,\nabla_q]\,\beta_{ij}=0
$$
and substitute the curvature tensor components \thetag{14.1} into this 
identity:
$$
4\,\tilde\beta_{kq}\,\beta_{ij}+\beta_{ki}\,\beta_{qj}-\beta_{qi}\,
\beta_{kj}+\beta_{kj}\,\beta_{iq}-\beta_{qj}\,\beta_{ik}=0.
$$
Hence $2\,\tilde\beta_{kq}\,\beta_{ij}+\tilde\beta_{ki}\,\beta_{qj}
-\tilde\beta_{qi}\,\beta_{kj}=0$. This relationship is skew-symmetric
respective to $k$ and $q$. We symmetrize it respective to $i$ and $j$:
$$
4\,\tilde\beta_{kq}\,\hat\beta_{ij}+\tilde\beta_{ki}\,\hat\beta_{qj}
+\tilde\beta_{kj}\,\hat\beta_{qi}-\tilde\beta_{qi}\,\hat\beta_{kj}
-\tilde\beta_{ki}\,\hat\beta_{qj}=0.\hskip -3em
\tag17.20
$$
In the above projectively-euclidean coordinates $y^1,\,\ldots,\,y^n$
the only nonzero components of tensor $\tilde{\boldsymbol\beta}$
for $m=2$ are $\tilde\beta_{12}$ and $\tilde\beta_{21}=-\tilde\beta_{12}$.
Tensor $\hat{\boldsymbol\beta}$ in these coordinates can have only four 
nonzero components, they are $\hat\beta_{11}$, $\hat\beta_{22}$, $\hat
\beta_{12}$, and $\hat\beta_{21}=\hat\beta_{12}$. For $k=1$ and $q=2$
from \thetag{17.20} we get
$$
\xalignat 3
&\quad\tilde\beta_{12}\,\hat\beta_{11}=0,&&\tilde\beta_{12}\,\hat
\beta_{12}=0,&&\tilde\beta_{12}\,\hat\beta_{22}=0.\hskip -3em
\tag17.21
\endxalignat
$$
Following three components $\hat\beta_{11}$, $\hat\beta_{12}$, and
$\hat\beta_{22}$ cannot vanish simultaneously. Therefore from
\thetag{17.21} we get $\tilde\beta_{12}=-\tilde\beta_{21}=0$, which
leads to $(n+1)\,\tilde{\boldsymbol\beta}=-\tilde\bold R=0$.
\proclaim{Lemma 17.1} Proposition of \pagebreak lemma~9.2 remains true 
for $m=2$.
\endproclaim
     Let's alternate \thetag{14.3} respective to $i$ and $j$ and let's
take into account symmetry of connection components. Due to $\beta_{ij}=
\beta_{ji}$ we obtain
$$
\frac{\partial u_j}{\partial y^i}-\frac{\partial u_i}{\partial y^j}=0.
\hskip -3em
\tag17.22
$$
The relationships \thetag{17.22} are exactly the compatibility conditions 
for the following system of Pfaff equations for the scalar field $\psi$:
$$
\frac{\partial\psi}{\partial y^i}=u_i(y^1,\ldots,y^m),\quad
i=1,\,\ldots,\,m.\hskip -3em
\tag17.23
$$
Due to the integrability of equations \thetag{17.23} we can state a
theorem which makes \thetag{17.19} more specific.
\proclaim{Theorem 17.3} Let \thetag{1.1} be a system of equations belonging 
to $m$-th case of intermediate degeneration such that its algebra of point 
symmetries has maximal dimension \thetag{9.1}. Under these assumptions 
there is a point transformation \thetag{1.4} bringing these equations to the 
form
$$
\frac{\partial y^i}{\partial\,t}=a\,\frac{\partial^2 y^i}{\partial x^2}
-2\,a\left(\,\shave{\sum^m_{r=1}}\frac{\partial\psi}{\partial y^r}
\,\frac{\partial y^r}{\partial x}\right)\frac{\partial y^i}{\partial x},
\quad i=1,\,\ldots,\,n,\hskip -3em
\tag17.24
$$
where $a=\const$ and $\psi=\psi(y^1,\ldots,y^m)$ is some scalar field.
\endproclaim
\head
18. Case of general position.
\endhead
     Case of general position is distinguished by the condition
$m=n\geqslant 2$, which means that the kernel of bilinear form
defined by symmetric part of Ricci tensor is zero: $\Ker\hat R=0$.
The estimate \thetag{8.9} for this case yields
$$
\dim(\Cal G)\leqslant\frac{n(n+1)}{2}.\hskip -3em
\tag18.1
$$
Our further aim is to describe systems of equations \thetag{1.1}
for which the upper bound in the estimate \thetag{18.1} is reached:
$$
\dim(\Cal G)=\frac{n(n+1)}{2}.\hskip -3em
\tag18.2
$$
Let's start with study of the fields $\bold S$ and $\tilde\bold R$ from 
\thetag{4.11} and \thetag{4.12}. For these fields we have a lemma
similar to lemma~9.2.
\proclaim{Lemma 18.1} Let \thetag{1.1} be a system of equations belonging 
to the case of general position such that its algebra of point symmetries 
has maximal dimension \thetag{18.2}. Under these assumptions $\bold S=2\,
\tilde\bold R=0$.
\endproclaim
     Proof of this lemma is almost literally the same as the proof of
lemma~9.2. The difference is only that $\Ker\hat R=0$. Therefore we
need not use factorspace and factoroperators in proving lemma~18.1.
\proclaim{Lemma 18.2} Let \thetag{1.1} be a system of equations belonging 
to the case of general position such that its algebra of point symmetries 
has maximal dimension \thetag{18.2}. Under these assumptions if $n\geqslant 
3$, then $\bold Q=\nabla\bold R=0$.
\endproclaim
\proclaim{Lemma 18.3} Proposition of lemma~18.2 remains true for the 
dimension $n=2$.
\endproclaim
     Proof of the lemma~18.2 is based on lemma~18.1. It is the same as
the proof of lemma~15.3 with the only difference $W=\Ker\hat R=0$.
Therefore we can do without factorspace $V/W$. Proof of lemma~18.3 is
analogous to that of the lemma~15.4.\par
     In the case of general position tensor $\hat\bold R$ is nondegenerate.
Therefore it defines pseudoriemannian metric $\bold g=\hat\bold R$
on $M$. From lemmas~18.2 and 18.3 we find that this metric is in
concordance with the connection $\Gamma$:
$$
\nabla_k g_{ij}=\nabla_k\hat R_{ij}=\frac{\nabla_kR_{ij}+
\nabla_kR_{ji}}{2}=0.\hskip -3em
\tag18.3
$$
Components of symmetric connection $\Gamma$ concordant with metric
tensor $\bold g$ are calculated by well-known formula 
$$
\Gamma^k_{ij}=\sum^n_{i=1}\frac{g^{ks}}{2}\left(\frac{\partial g_{sj}}
{\partial y^i}+\frac{\partial g_{is}}{\partial y^j}-\frac{\partial g_{ij}}
{\partial y^s}\right).\hskip -3em
\tag18.4
$$
Formula \thetag{18.4} for $\Gamma^k_{ij}$ follows from \thetag{18.3}
(details see in \cite{6}, \cite{9}, or \cite{17}).
     Using metric tensor $\bold g$ one can lower upper index in curvature
tensor. This leads to the tensor field of type $(0,4)$ with components
$$
R_{kqij}=\sum^n_{i=1}g_{ks}\,R^s_{qij}.\hskip -3em
\tag18.5
$$
\proclaim{Theorem 18.1} Components of curvature tensor \thetag{18.5}
for metric connection \thetag{18.4} possess two symmetry properties
$$
\xalignat 2
&R_{kqij}=-R_{qkij},&&R_{kqij}=R_{ijkq},\hskip -3em
\tag18.6
\endxalignat
$$
in addition to the properties $R^k_{qij}=-R^k_{qji}$ and $R^k_{qij}
+R^k_{ijq}+R^k_{jqi}=0$, which are available for curvature tensor of
any symmetric connection.
\endproclaim
     Proposition of theorem~18.1 is well-known fact. Its proof can be
found in \cite{6} or in \cite{17}. As an immediate consequence of 
\thetag{18.6} we obtain symmetry of Ricci tensor:
$$
R_{qj}=\sum^n_{i=1}\sum^n_{k=1}g^{ki}\,R_{kqij}=\sum^n_{i=1}\sum^n_{k=1}
g^{ik}\,R_{ijkq}=R_{jq}.
$$
This relationship can be used to strengthen the proposition of lemma~18.1
by extending it for the case of two-dimensional manifolds.
\proclaim{Lemma 18.4} Proposition of lemma~18.1 remains true for the 
dimension $n=2$.
\endproclaim
\head
19. Structure of curvature tensor for $n\geqslant 3$.
\endhead
    Let's fix a point $p_0$ on the manifold $M$ and consider a field of 
point symmetry $\boldsymbol\eta$ of the system of equations \thetag{1.1} 
vanishing at the point $p_0$. Such field generates linear operator
\thetag{5.7} in tangent space $T_{p_0}(M)$. Denote it by $\bold F$. 
According to the theorem~4.1 we have $L_{\boldsymbol\eta}(\bold R)=0$.
This equation written at the point $p_0$ takes the form \thetag{10.2}.
Here we suppose that system of equations \thetag{1.1} belongs to the 
case of general position and the dimension of its symmetry algebra 
is given by \thetag{18.2}. Due to \thetag{18.2} any operator $\bold F$ 
satisfying \thetag{8.3} corresponds to some symmetry field $\boldsymbol
\eta\in\Cal G(p_0)$. In the case of general position bilinear form 
$g=\hat R$ in \thetag{8.3} is nondegenerate. By applying the 
complexification procedure to the space $V$, if necessary, the matrix of 
the form $g$ can be brought to the unit matrix at the expense of proper 
choice of base $\bold e_1,\,\ldots,\,\bold e_n$ in $V=T_{p_0}(M)$:
$$
g_{ij}=g(\bold e_i,\bold e_j)=\hat R(\bold e_i,\bold e_j)=\cases 1
&\text{for \ }i=j,\\ 0 &\text{for \ }i\neq j.\endcases\hskip -3em
\tag19.1
$$
For the dimension $n\geqslant 3$ we can consider three base vectors
$\bold e_k$, $\bold e_q$, $\bold e_r$ such that $k\neq q\neq r\neq k$.
Let's define an operator $\bold F$ by its action on the base vectors:
$$
\hat\bold F(\bold e_i)=\cases\hphantom{-}\bold e_q &\text{for \ }i=k,\\
-\bold e_k &\text{for \ }i=q,\\ \hphantom{-}0 &\text{for \ }i\neq k
\text{\ \ and \ }i\neq q.\endcases\hskip -3em
\tag19.2
$$
One can easily check that for operator $\bold F$ the relationships 
\thetag{8.3} are fulfilled. Let's substitute $\bold F$ into the
formula \thetag{10.2} and take $\bold X=\bold e_k$, $\bold Y=\bold e_r$, 
$\bold Z=\bold e_k$ in it. This yields $\bold R(\bold F\bold e_k,
\bold e_r)\bold e_k+\bold R(\bold e_k,\bold e_r)\bold F\bold e_k
=\bold F\bold R(\bold e_k,\bold e_r)\bold e_k$. Now by \thetag{19.2}
we obtain
$$
\sum^n_{i=1}R^i_{kqr}\,\bold e_i+\sum^n_{i=1}R^i_{qkr}\,\bold e_i=
\sum^n_{s=1}R^s_{kkr}\,\bold F\bold e_s=R^k_{kkr}\,\bold e_q-
R^q_{kkr}\,\bold e_k.
$$
This vectorial equality is reduced to the series of scalar equalities:
$$
\align
&R^k_{kqr}+R^k_{qkr}=-R^q_{kkr}\text{\ \ for \ }i=k,\hskip -3em
\tag19.3\\
\vspace{1ex}
&R^q_{kqr}+R^q_{qkr}=R^k_{kkr}\text{\ \ for \ }i=q,\hskip -3em
\tag19.4\\
\vspace{1ex}
&R^i_{kqr}+R^i_{qkr}=0\text{\ \ for \ }i\neq k\text{\ \ and \ }i\neq q.
\hskip -3em
\tag19.5
\endalign
$$
Due to \thetag{19.1} we make no difference between upper and lower 
indices for tensor components referred to the base $\bold e_1,\,\ldots,
\,\bold e_n$. By means of \thetag{18.6} from \thetag{19.3}, 
\thetag{19.4}, and \thetag{19.5} we derive the following 
equalities:
$$
R^k_{kqr}=0,\quad R^q_{kqr}=0\text{, \ and \ }R^i_{kqr}=-R^i_{qkr}
\text{\ \ for \ }i\neq k\text{\ \ and \ }i\neq q.
$$
Note that we can unite these three equalities into one:
$$
R^i_{kqr}=-R^i_{qkr}\text{\ \ for \ }k\neq q\neq r\neq k.\hskip -3em
\tag19.6
$$
Let's take into account \thetag{19.6} and write the identity $R^i_{kqr}
+R^i_{qrk}+R^i_{rkq}=0$, which follows from symmetry $\Gamma^i_{rs}=
\Gamma^i_{sr}$. This results in
$$
R^i_{kqr}=0\text{\ \ for \ }q\neq k\neq r.\hskip -3em
\tag19.7
$$\par
     As a next step let's substitute $\bold F$ into the formula
\thetag{10.2} taking $\bold X=\bold e_k$, $\bold Y=\bold e_r$, 
$\bold Z=\bold e_q$ in it: $\bold R(\bold F\bold e_k,\bold e_r)
\bold e_q+\bold R(\bold e_k,\bold e_r)\bold F\bold e_q=\bold F
\bold R(\bold e_k,\bold e_r)\bold e_q$. By \thetag{19.2} we get
$$
\sum^n_{i=1}R^i_{qqr}\,\bold e_i-\sum^n_{i=1}R^i_{kkr}\,\bold e_i=
\sum^n_{s=1}R^s_{qkr}\,\bold F\bold e_s=R^k_{qkr}\,\bold e_q-
R^q_{qkr}\,\bold e_k.
$$
This vectorial equality is reduced to the series of scalar equalities:
$$
\align
&R^k_{qqr}-R^k_{kkr}=-R^q_{qkr}\text{\ \ for \ }i=k,\\
\vspace{1ex}
&R^q_{qqr}-R^q_{kkr}=R^k_{qkr}\text{\ \ for \ }i=q,\\
\vspace{1ex}
&R^i_{qqr}-R^i_{kkr}=0\text{\ \ for \ }i\neq k\text{\ \ and \ }i\neq q.
\endalign
$$
We choose the last equality and substitute $i=r$ in it. As a result
we obtain
$$
R^r_{qqr}=R^r_{kkr}\text{\ \ for \ }q\neq r\neq k.\hskip -3em
\tag19.8
$$\par
     Due to \thetag{19.7} we conclude that $R^k_{qsr}$ can be nonzero
only if some two of three lower indices are equal, i\.~e\. if $q=r$ 
or $q=s$ (chance $s=r$ is excluded since $R^k_{qsr}$ is skew-symmetric
relative to indices $s$ and $r$):
$$
R^k_{qsr}=g_{qr}\,\beta^k_s(q)-g_{qs}\,\beta^k_r(q).\hskip -3em
\tag19.9
$$
In formula \thetag{19.9} we took into account \thetag{19.1} and
skew-symmetry $R^k_{qsr}=-R^k_{qrs}$. Now we are to specify the
quantities $\beta^k_s(q)$ and $\beta^k_r(q)$ in \thetag{19.9}.\par
     Let $q\neq r$. In order to find $\beta^k_r(q)$ we substitute 
$s=q$ into the formula \thetag{19.9}. Upon some easy calculations
we get
$$
\beta^k_r(q)=-R^k_{qqr}=R_{kqqr}=R_{qkqr}=R^q_{kqr}\text{\ \ for \ }
q\neq r.\hskip -3em
\tag19.10
$$
Here we took into account \thetag{18.6} and the property \thetag{19.1}
of the base $\bold e_1,\,\ldots,\,\bold e_n$, due to which we can mix
upper and lower indices. Let's transform the right hand side of
\thetag{19.10}. For $q\neq k\neq r$ right hand side of \thetag{19.10}
is zero, one should substitute $i=q$ into \thetag{19.7} in order
to see it. For $q=k\neq r$ it is zero too due to $R^q_{qqr}=R_{qqqr}$
and \thetag{18.6}. Therefore
$$
\beta^k_r(q)=0\text{\ \ for \ }q\neq r\neq k.\hskip -3em
\tag19.11
$$
Suppose $q\neq k=r$. Then from \thetag{19.10} we get $\beta^r_r(q)=
-R^r_{qqr}$. \pagebreak Formula \thetag{19.8} and the equalities 
$R^r_{rrr}=0$ and $R^r_{kkr}=R_{rkkr}=-R_{krkr}=-R^k_{rkr}$ are used 
for further transformation of the expression for $\beta^r_r(q)$:
$$
\beta^r_r(q)=-R^r_{qqr}=-\sum^n_{k\neq r}\frac{R^r_{kkr}}{n-1}=
\sum^n_{k=1}\frac{R^k_{rkr}}{n-1}=\frac{R_{rr}}{n-1}=\frac{g_{rr}}
{n-1}=\frac{1}{n-1}.
$$
We combine this result with \thetag{19.11} and we write it as follows:
$$
\beta^k_r(q)=\frac{\delta^k_r}{n-1}\text{\ \ for \ }q\neq r.\hskip -3em
\tag19.12
$$
Formula \thetag{19.12} do not define $\beta^k_r(q)$ for $q=r$. For 
arbitrary values of $k$, $q$, and $r$ we replace it by the following 
formula, equivalent to \thetag{19.12} for $q\neq r$:
$$
\beta^k_r(q)=\frac{\delta^k_r}{n-1}(1-g_{qr})+\beta^k_r(r)\,g_{qr}.
\hskip -3em
\tag19.13
$$
Now substitute \thetag{19.13} into the formula \thetag{19.9}. As a result
we get
$$
R^k_{qsr}=\frac{g_{qr}\,\delta^k_s-g_{qs}\,\delta^k_r}{n-1}.\hskip -3em
\tag19.14
$$
It's worth to compare \thetag{19.14} with \thetag{10.15} and take into
account lemma~18.1, due to which Ricci tensor $\bold g=\bold R$ is 
symmetric.
\head
20. Structure of curvature tensor for $n=2$.
\endhead
     For $n=2$ Ricci tensor is symmetric as well. This follows from
lemma~18.4. Formula \thetag{19.14} for $n=2$ takes the following form:
$$
R^k_{qsr}=g_{qr}\,\delta^k_s-g_{qs}\,\delta^k_r.\hskip -3em
\tag20.1
$$
The proof or the formula \thetag{20.1} is the same as fore the formula
\thetag{12.1}.
\proclaim{Theorem 20.1} Let \thetag{1.1} be a system of equations belonging 
to the case of general position such that its algebra of point symmetries 
has maximal dimension \thetag{18.2}. Under these assumptions if $n=2$,
then components of curvature tensor are expressed through components 
Ricci tensor $\bold g=\bold R$ according to the formula \thetag{19.4}.
\endproclaim
\head
21. Structure of tensor field $\bold A$
\endhead
\proclaim{Theorem 21.1} Let \thetag{1.1} be a system of equations belonging 
to the case of general position such that its algebra of point symmetries 
has maximal dimension \thetag{18.2}. Under these assumptions if $n\geqslant 
3$, then $A$ is a constant scalar matrix, i\.~e\. $A^i_j=a\,\delta^i_j$
and $a=\const$.
\endproclaim
\demo{Proof} Let's fix a point $p_0$ on $M$ \pagebreak and consider a field 
of point symmetry $\boldsymbol\eta$ for the system of equations \thetag{1.1}. 
Suppose that $\boldsymbol\eta$ vanishes at $p_0$. The equation
$L_{\boldsymbol\eta}(\bold A)=0$ from \thetag{2.5} for this field at the
point $p_0$ is written as
$$
\bold A\,\bold F\bold X=\bold F\,\bold A\bold X.\hskip -3em
\tag21.1
$$
Here $\bold F$ is linear operator from \thetag{5.7} and $\bold X$
is an arbitrary vector from tangent space $V=T_{p_0}(M)$. Due to
\thetag{18.2} for $\bold F$ we can take any operator satisfying the 
equations \thetag{8.3}. In the case of general position bilinear form 
$g=\hat R$ in \thetag{8.3} is nondegenerate. By applying the 
complexification procedure to the space $V$, if necessary, the matrix 
of nondegenerate symmetric bilinear form $\hat R$ from \thetag{8.7} 
can be brought to the unit matrix at the expense of proper choice of 
base $\bold e_1,\,\ldots,\,\bold e_m$, i\.~e\. there is a base such
that the relationships \thetag{19.1} hold.\par
     For the dimension $n\geqslant 3$ we can find three base vectors
$\bold e_k$, $\bold e_q$, $\bold e_r$ such that $k\neq q\neq r\neq k$.
Let's define the operator $\bold F$ by means of \thetag{19.2} and 
substitute it into the relationship \thetag{21.1} taking $\bold X
=\bold e_k$. This yields
$$
\sum^n_{i=1}A^i_q\,\bold e_i=\sum^n_{s=1}A^s_k\,\bold F\bold e_s=
A^k_k\,\bold e_q-A^q_k\,\bold e_k.
$$
The above vectorial equality reduces to the series if scalar equalities:
$$
A^q_q=A^k_k,\quad A^k_q=-A^q_k\text{, \ and \ }A^i_q=0\text{\ \ for \ }
i\neq k\text{\ \ and \ }i\neq q.\hskip -3em
\tag21.2
$$
First equality $A^q_q=A^k_k$ means, that all diagonal elements of 
matrix $A$ are equal to each other. The last equality $A^i_q=0$ means
that nondiagonal elements are zero. Hence $A^i_j=a\,\delta^i_j$ or
$\bold A=a\,\id$, where $a$ is a scalar field. Substituting $\bold A=a\,
\id$ into the equation $L_{\boldsymbol\eta}(\bold A)=0$ we get the
equation \thetag{13.4}, which leads to $L_{\boldsymbol\eta}(a)=0$.
Due to \thetag{18.2} we can find $n$ vector fields $\boldsymbol\eta_1,
\,\ldots,\,\boldsymbol\eta_n$ of point symmetry algebra $\Cal G$
which are linear independent at the point $p_0$. The equation
$L_{\boldsymbol\eta}(a)=0$ holds for each of these fields. Therefore
$a=\const$. Theorem is proved.\qed\enddemo
     On the orientable manifold of the dimension $n=2$ equipped with
nondegenerate metric $\bold g$ apart from the field of identical
operators $\id$ there exists a field of rotation by the angle
$90^\circ$ in metric $\bold g$. Denote it by $\bold P$. Components
of tensor field $\bold P$ are defined by the following explicit formula:
$$
P^i_j=\pm\sum^2_{s=1}\frac{d^{is}\,g_{sj}}{\sqrt{|\det g|}}.\hskip -3em
\tag21.3
$$
Here $d^{is}$ are components of skew-symmetric unit matrix \thetag{12.2}.
Sign in the formula \thetag{21.3} is defined by the orientation of the
system of local coordinates. Upon complexification of the tangent space
$V$ if we choose the base where the relationships \thetag{19.1} are
fulfilled, then for the matrix of the operator field $\bold P$ we obtain
$$
\pagebreak
P=\varepsilon\cdot\Vmatrix 0 & 1\\ \vspace{1ex}-1 & 0 \endVmatrix,
\hskip -3em
\tag21.4
$$
where $\varepsilon=\pm 1$ for positive or negative metric $\bold g$
and $\varepsilon=\pm i$ for indefinite metric $\bold g$.
\proclaim{Theorem 21.2} Let \thetag{1.1} be a system of equations belonging 
to the case of general position such that its algebra of point symmetries 
has maximal dimension \thetag{18.2}. If $n=2$, then components of the matrix 
$A$ are given by the formula
$$
A^i_j=a\,\delta^i_j+b\,\sum^2_{s=1}\frac{d^{is}\,R_{sj}}
{\sqrt{|\det R|}},\hskip -3em
\tag21.5
$$
where $a$ and $b$ are constants, $R_{sj}=g_{sj}$ are components of Ricci 
tensor for the connection $\Gamma$, and $\det R$ is the determinant of 
the matrix of this tensor.
\endproclaim
\demo{Proof} For the dimension $n=2$ we have only two base vectors 
$\bold e_1$ and $\bold e_2$. Let's define the operator $\bold F$ by
its action on base vectors:
$$
\bold F(\bold e_i)=\cases\hphantom{-}\bold e_2 &\text{for \ }
i=1,\\ -\bold e_1 &\text{for \ }i=2.\endcases
$$
When we substitute $\bold F$ into \thetag{21.1} and take $\bold X=\bold 
e_1$, instead of \thetag{21.2} here we get $A^2_2=A^1_1=a$ and $A^2_1=
-A^1_2=-b$. Therefore matrix of operator $\bold A$ has the form
$$
A=\Vmatrix a & b\\ \vspace{1ex}-b & a\endVmatrix.\hskip -3em
\tag21.6
$$
Comparing \thetag{21.4} with \thetag{21.6} we obtain \thetag{21.5}.
Furthermore, from \thetag{21.6} we get $2\,a=\tr\bold A$ and $a^2
+b^2=\det\bold A$. Therefore the equation $L_{\boldsymbol\eta}(
\bold A)=0$ yields $L_{\boldsymbol\eta}(a)=0$ and $L_{\boldsymbol
\eta}(b)=0$. Due to \thetag{18.2} then we have $a=\const$ and 
$b=\const$.\qed\enddemo
\head
22. Constant curvature spaces.
\endhead
     If the system of equation \thetag{1.1} belongs to the case of
general position and if its point symmetry algebra has maximal
dimension \thetag{18.2}, then curvature tensor for appropriate
affine connection $\Gamma$ has special structure \thetag{19.14},
which corresponds to the spaces of constant sectional curvature
$K=1/(n-1)$ (see \cite{24}). In special local coordinates metric
tensor and components of connection for such spaces can be brought
the special form. In order to construct such coordinates let's
consider the following system of Pfaff equations for the components
of covector field $\bold u$:
$$
\gather
\nabla_iu_j=-u_i\,u_j-\frac{g_{ij}}{2(n-1)}+\frac{|\bold u|^2}{2}
\,g_{ij},\quad i,j=1,\,\ldots,\,n,\hskip -3em
\tag22.1\\
\vspace{1ex}
\text{where \ }|\bold u|^2=\sum^n_{r=1}\sum^n_{s=1}g^{rs}\,u_r\,u_s.
\hskip -3em
\endgather
$$
Complete compatibility of Pfaff equations \thetag{22.1} follows from 
\thetag{19.14} and \thetag{18.3}. For the equations \thetag{22.1}
let's consider Cauchy problem with zero initial data
$$
u_j\,\hbox{\vrule height 6pt depth 8pt width 0.5pt}_{\,p=p_0}
=0\hskip -3em
\tag22.2
$$
at some fixed point $p_0$. Right hand side of the equations \thetag{22.1}
do not vanish for $\bold u=0$. Therefore the solution of the Cauchy
problem \thetag{22.2} is unique covector field, which vanishes at the
point $p_0$, but which is not identically zero.\par
     By raising index with the use of metric $\bold g$ we convert
covector field $\bold u$ into the vector field $\bold u$ with
components
$$
u^k=\sum^n_{j=1}g^{kj}\,u_j.
$$
It satisfies the system of Pfaff equations derived from \thetag{22.1}:
$$
\nabla_iu^j=-u_i\,u^j-\frac{\delta^j_i}{2(n-1)}+\frac{|\bold u|^2}{2}
\,\delta^j_i,\quad i,j=1,\,\ldots,\,n.\hskip -3em
\tag22.3
$$
By means of components of the field $\bold u$ we construct a tensor
field $\bold T$ of type $(1,2)$. Its components are defined as follows:
$$
T^k_{rs}=-u_r\,\delta^k_s-u_s\,\delta^k_r+u^k\,g_{rs}.\hskip -3em
\tag22.4
$$
We use tensor $\bold T$ as a deformation tensor for connection $\Gamma$:
$$
\bar\Gamma^k_{rs}=\Gamma^k_{rs}+T^k_{rs}=\Gamma^k_{rs}
-u_r\,\delta^k_s-u_s\,\delta^k_r+u^k\,g_{rs}.\hskip -3em
\tag22.5
$$
Its not difficult to calculate curvature tensor for new connection
$\bar\Gamma$:
$$
\bar R^k_{sij}=R^k_{sij}+\nabla_iT^k_{js}-\nabla_jT^k_{is}+\sum^n_{q=1}
T^q_{js}\,T^k_{iq}-\sum^n_{q=1}T^q_{js}\,T^k_{iq}.\hskip -3em
\tag22.6
$$
Substituting \thetag{19.14} and \thetag{22.4} into \thetag{22.6}
and taking into account the equations \thetag{22.1} and \thetag{22.3} 
we get $\bar R^k_{sij}=0$. Thus \thetag{22.5} is a plane (euclidean)
connection.\par
     Right hand side of the equations \thetag{22.1} is symmetric
respective to $i$ and $j$. Therefore $\nabla_iu_j=\nabla_ju_i$. Due
to the symmetry $\Gamma^k_{ij}=\Gamma^k_{ji}$ this reduces to
$$
\frac{\partial u_j}{\partial y^i}=\frac{\partial u_i}{\partial y^j}.
$$
These relationships mean that covector field $\bold u$ is a gradient
of some scalar field. One can find an explicit formula for this
field. Let's take 
$$
f=\frac{1}{2(n-1)}+\frac{|\bold u|^2}{2}.\hskip -3em
\tag22.7
$$
Function $f$ in \thetag{22.7} is positive in some neighborhood of
the point $p_0$ which defines initial data for Cauchy problem 
\thetag{22.2}. From \thetag{22.1} we derive
$$
\nabla_if=-u_i\,f.\hskip -3em
\tag22.8
$$
So covector field $\bold u$ \pagebreak is a gradient of the scalar field 
$\psi=-\ln f$.\par
     Let's consider the metric $\bar\bold g=f^2\,\bold g$ conformally
equivalent to the initial metric $\bold g=\bold R$. By means of direct
calculations with the use of formula  \thetag{18.4}  we  can  see 
that
metric connection for new metric $\bar\bold g=f^2\,\bold g$ coincides
with  \thetag{22.5}.
\proclaim{Conclusion} Metric $\bar\bold g=f^2\,\bold g$ is a flat
(pseudoeuclidean) metric and there exist some local coordinates $y^1,\,
\ldots,\,y^n$ on $M$ such that components of the metric $\bar\bold g$ 
are constants and components of the connection $\bar\Gamma$ are zero.
\endproclaim
     With respect to initial metric $\bold g$ such coordinates are called
{\it conformally-euclidean} coordinates. In conformally-euclidean
coordinates we have
$$
\Gamma^k_{rs}=u_r\,\delta^k_s+u_s\,\delta^k_r-u^k\,g_{rs}.\hskip -3em
\tag22.9
$$
These coordinates can be chosen so that constant matrix $\bar g_{ij}$
is brought to the canonical form. Then for initial metric $\bold g$
we get
$$
\bold g=\sum^n_{i=1}\,\frac{\varepsilon_i\,dy^i\otimes dy^i}{f^2},
\hskip -3em
\tag22.10
$$
where $\varepsilon_1=\pm 1,\,\ldots,\,\varepsilon_n=\pm 1$. The number
of pluses and minuses in the sequence $\varepsilon_1,\,\ldots,\,
\varepsilon_n$ is defined by the signature of Ricci tensor $\bold R=
\bold g$.\par
     In order to find explicit form for the metric \thetag{22.10} we
are only to find the explicit form of the function $f$ in
conformally-euclidean coordinates. Let's write \thetag{22.3} and 
\thetag{22.8} in these coordinates and take into account \thetag{22.7} 
and \thetag{22.9}:
$$
\xalignat 2
&\frac{\partial f}{\partial y^i}=-u_i\,f,
&&\frac{\partial u^j}{\partial y^i}=-u_i\,u^j-f\,\delta^j_i.
\hskip -3em
\tag22.11
\endxalignat
$$
Now let's resolve the first equation \thetag{22.11} with respect to 
$u_i$ and substitute the obtained expression for $u_i$ into the second
equation \thetag{2.11}. This yields
$$
\frac{\partial}{\partial y^i}\left(-\frac{u^j}{f}\right)=\delta^j_i.
$$
This is the system of Pfaff equations which can be integrated  in 
explicit
form. Since initial data in \thetag{22.2} are zero, we get $u^j=-f\,y^j$. 
Here we assumed that fixed point $p_0$ is the origin for 
conformally-euclidean coordinates $y^1,\,\ldots,\,y^n$. Let's substitute
$u^j=-f\,y^j$ into \thetag{22.7} and take into account formula 
\thetag{22.10} for metric. As a result for $f$ we obtain explicit formula
in conformally-euclidean coordinates:
$$
f=\frac{1}{2(n-1)}+\frac{1}{2}\sum^n_{i=1}\varepsilon_i\,(y^i)^2.
\hskip -3em
\tag22.12
$$
So we got the proof of the following theorem first proved by Riemann.
\proclaim{Theorem 22.1} In some neighborhood of any point $p_0$
on $n$-dimensional pseudoriemannian manifold of constant sectional 
curvature \pagebreak $K=1/(n-1)$ there exist conformally-euclidean 
coordinates $y^1,\,\ldots,\,y^n$ such that metric tensor has the form 
\thetag{22.10} with the parameter $f$ defined by \thetag{22.12}.
\endproclaim
\proclaim{Theorem 22.2} Let \thetag{1.1} be a system of equations belonging 
to the case of general position such that its algebra of point symmetries 
has maximal dimension \thetag{18.2}. If $n\geqslant 3$, there is a point 
transformation \thetag{1.4} bringing these equations to the form
$$
\frac{\partial y^i}{\partial\,t}=a\,\frac{\partial^2 y^i}{\partial x^2}
+\sum^n_{r=1}a\,\frac{\varepsilon_r}{f}\,\frac{\partial y^r}{\partial x}
\left(y^i\,\frac{\partial y^r}{\partial x}-2\,y^r\frac{\partial y^i}
{\partial x}\right),\hskip -3em
\tag22.13
$$
where $i=1,\,\ldots,\,n$, $a=\const$, $\varepsilon_1=\pm 1,\,\ldots,\,
\varepsilon_n=\pm 1$, and $f=f(y^1,\ldots,y^n)$ is a function defined
by \thetag{22.12}.
\endproclaim
    The equations \thetag{22.13} arise if we substitute \thetag{22.9}
into \thetag{1.1} and take into account $u^j=-f\,y^j$. In two-dimensional
case $n=2$ appropriate equations appear to be more huge due to more
huge formula for $\bold A$:
$$
\align
&\aligned
\frac{\partial y^1}{\partial\,t}&=a\,\frac{\partial^2 y^1}{\partial x^2}
+\sum^n_{r=1}a\,\frac{\varepsilon_r}{f}\,\frac{\partial y^r}{\partial x}
\left(y^1\,\frac{\partial y^r}{\partial x}-2\,y^r\frac{\partial y^1}
{\partial x}\right)\,+\\
&+\,b\,\varepsilon_1\,\frac{\partial^2 y^2}{\partial x^2}
+\sum^n_{r=1}b\,\varepsilon_1\,\frac{\varepsilon_r}{f}\,\frac{\partial y^r}
{\partial x}\left(y^2\,\frac{\partial y^r}{\partial x}-2\,y^r\frac{\partial
y^2}{\partial x}\right),
\endaligned\hskip -3em
\tag22.14\\
\vspace{2ex}
&\aligned
\frac{\partial y^2}{\partial\,t}&=a\,\frac{\partial^2 y^2}{\partial x^2}
+\sum^n_{r=1}a\,\frac{\varepsilon_r}{f}\,\frac{\partial y^r}{\partial x}
\left(y^2\,\frac{\partial y^r}{\partial x}-2\,y^r\frac{\partial y^2}
{\partial x}\right)-\\
&-b\,\varepsilon_2\,\frac{\partial^2 y^1}{\partial x^2}
-\sum^n_{r=1}b\,\varepsilon_2\,\frac{\varepsilon_r}{f}\,\frac{\partial y^r}
{\partial x}\left(y^1\,\frac{\partial y^r}{\partial x}-2\,y^r\frac{\partial
y^1}{\partial x}\right).
\endaligned\hskip -3em
\tag22.15
\endalign
$$
\proclaim{Theorem 22.3} Let \thetag{1.1} be a system of equations belonging 
to the case of general position such that its algebra of point symmetries 
has maximal dimension \thetag{18.2}. If $n=2$, then there is a point 
transformation \thetag{1.4} bringing these equations to the form
\thetag{22.14} and \thetag{22.15}, where $a=\const$, $b=\const$, 
$\varepsilon_1=\pm 1$, $\varepsilon_2=\pm 1$, and $f=f(y^1,y^2)$ is
a function defined by \thetag{22.12}.
\endproclaim
\head
23. Case of one equation.
\endhead
     For $n=1$ we have only one equation in the system \thetag{1.1}.
Operator field $\bold A$ and affine connection $\Gamma$ have only
one component:
$$
\frac{\partial y}{\partial\,t}=A\left(\frac{\partial^2 y}{\partial x^2}
+\Gamma\,\frac{\partial y}{\partial x}\,\frac{\partial y}{\partial x}
\right).\hskip -3em
\tag23.1
$$
Any affine connection on one-dimensional manifold has identically
zero curvature tensor. Therefore it is flat and its Ricci tensor
is zero. This means $m=0$. Any equation \thetag{23.1} belongs to
the case of maximal degeneration. For the dimension of its symmetry
algebra we have an estimate $\dim(\Cal G)\leqslant 2$, which follows
from theorem~6.1. 
\proclaim{Theorem 23.1} For the equation \thetag{23.1} with the
symmetry algebra of maximal dimension $\dim(\Cal G)=2$ there exist
a point transformation bringing it to the form
$$
\frac{\partial y}{\partial\,t}=a\,\frac{\partial^2 y}{\partial\,x^2}
\text{, \ \ where \ \ }a=\const.
$$
\endproclaim
{\bf Note}. According to theorem~17.2 if system of equations
\thetag{1.1} belongs to $m$-th case of intermediate degeneration
and possess the algebra of point symmetries of maximal dimension
\thetag{9.1}, then it admit the variable separation resulting in
smaller subsystem of $m$ equations. Such subsystem belongs to the
case of general position and has the algebra of point symmetries
of maximal dimension \thetag{18.2}, where $n=m$. To this system we
can apply one of the theorems~22.2, 22.3, or 23.1.\par

\enddocument
\end